\documentclass{article}
\usepackage{arxiv}
\usepackage{epsfig}
\usepackage{float}
\usepackage{amsmath,amssymb}
\usepackage{algorithm,algorithmic}
\usepackage[labelformat=simple]{subfig}

\usepackage{xcolor}
\usepackage{epstopdf}
\usepackage{multirow}
\usepackage{latexsym,bm}
\usepackage{doi}     % for things like \Box and \Diamond
\graphicspath{{./Figures/}}
\usepackage{hyperref}
\makeindex

\title{Energy preserving reduced-order modelling of thermal  shallow water equation}

\author{ S\"uleyman Y{\i}ld{\i}z\\
Institute of Applied Mathematics\\
Middle East Technical University\\
 Ankara-Turkey\\
	\texttt{yildiz.suleyman@metu.edu.tr} \\
	\And
  Murat Uzunca\\
   Department of Mathematics, Sinop University\\
     Sinop-Turkey\\ \
     \texttt{muzunca@sinop.edu.tr}\\
     \And
     B\"ulent Karas\"ozen\\
     Institute of Applied Mathematics \& Department of Mathematics\\
     Middle East Technical University\\
     Ankara-Turkey\\
     \texttt{bulent@metu.edu.tr}
}

\begin{document}
\maketitle

\date{Received: date / Accepted: date}

\begin{abstract}
 In this paper, Hamiltonian and energy preserving reduced-order models are developed for the rotating thermal shallow water equation (RTSWE) in the non-canonical Hamiltonian form with the state-dependent Poisson matrix. The high fidelity full solutions are obtained by discretizing the RTSWE in space with skew-symmetric finite-differences, that preserve the Hamiltonian structure. The resulting skew-gradient system is integrated in time with the energy preserving average vector field (AVF) method. The reduced-order model (ROM) is constructed in the same way as the full order model (FOM), preserving the reduced skew-symmetric structure and integrating in time with the AVF method. Relying on structure-preserving discretizations in space and time and applying proper orthogonal decomposition (POD) with the Galerkin projection, an energy preserving reduced order model (ROM) is constructed. The nonlinearities in the ROM are computed by applying the discrete empirical interpolation (DEIM) method to reduce the computational cost. The computation of the reduced-order solutions is accelerated further by the use of tensor techniques. The overall procedure yields a clear separation of the offline and online computational cost of the reduced solutions. The accuracy and computational efficiency of the ROMs are demonstrated for a numerical test problem. Preservation of the energy (Hamiltonian), and other conserved quantities, i.e. mass, buoyancy, and total vorticity show that the reduced-order solutions ensure the long-term stability of the solutions while exhibiting several orders of magnitude computational speedup over the FOM.
\end{abstract}

\keywords{Hamiltonian systems,   finite differences, proper orthogonal decomposition,  discrete empirical interpolation, tensors\\
MR2000 Subject Classification: 65M06; 65P10; 37J05; 37M15;  76B15}

%%%%%%%%%%%%%%%%%%%%%%%%%%%%%%%%%%%%%%%%%%%%%%%%%%%%%%%%%%%%%%%%%%%%%%%%%%%%%%%%%%%%%%%
%%%%%%%%%%%%%%%%%%%%%%%%%%%%%%%%%%%%%%%%%%%%%%%%%%%%%%%%%%%%%%%%%%%%%%%%%%%%%%%%%%%%%%%
%%%%%%%%%%%%%%%%%%%%%%%%%%%%%%%%%%%%%%%%%%%%%%%%%%%%%%%%%%%%%%%%%%%%%%%%%%%%%%%%%%%%%%%
\section{Introduction}

Rotating shallow water  equation (RSWE) \cite{Salmon04} is a  widely used  conceptual   model in geophysical and planetary fluid dynamics for the behavior of rotating inviscid fluids with one or more layers. Within each layer the horizontal velocity is assumed to be depth-independent, so the fluid moves in columns. However, the RSWE model does not
allow for gradients of the mean temperature and/or density, which are ubiquitous in the
atmosphere and oceans. The rotating thermal shallow water equation (RTSWE)  \cite{Dellar03,Eldred19,Ripa95,Dellar13} represents an extension of the RSWE, to include horizontal density/temperature gradients  both in the atmospheric and oceanic context. The RTSWE is used in the general circulation models \cite{Zerroukat15}, planetary flows \cite{Dellar14}, and modelling atmospheric and oceanic temperature fronts \cite{Dempsey88,Young95}, thermal instabilities \cite{Gouzien17}.

Shallow water equations (SWEs) are discretized on fine space-time grids to obtain high fidelity solutions. Real-time simulations require a large amount of computer memory and computing time. Therefore the computational cost associated with fully resolved simulations remains a barrier in many applications.
Model  order  reduction (MOR)  allows  to construct low-dimensional reduced-order models (ROMs) for the high dimensional full-order models (FOMs),
 generated by the discretization  of partial differential equations (PDEs) with  the finite difference, finite-volume, finite-element, spectral elements, discontinuous Galerkin methods.  The ROMs are computationally efficient and accurate and are worthy when a FOM needs to be simulated multiple-times for different parameter settings or in  multi-query scenarios such as  in  optimization. Additionally, ROMs are even more valuable for SWEs in simulating and predicting the model for a long time horizon.  The solutions of the high fidelity FOMs are projected on low dimensional reduced spaces usually  using the proper orthogonal decomposition (POD) \cite{Berkooz93,Sirovich87}, which is a widely used ROM technique. Applying the POD with the Galerkin projection, the dominant POD modes are extracted from the snapshots of the FOM solutions.
An important feature of the ROMs is the offline-online decomposition.  The computation of the FOM  and the construction of the reduced basis are  performed in the offline stage, whereas the reduced system is solved by projecting the problem onto   the low-dimensional reduced space in the online stage.  We refer to the books \cite{Benner17,Rozza14} for an overview of the available techniques. MOR techniques  for SWEs have been intensively studied in the literature, see, e.g., \cite{Navon20,Bistrian15,Bistrian17,Esfahanian09,Karasozen21,Lozovskiy17,Lozovskiy16,cstefuanescu2014comparison}.

Conservation of nonlinear invariants  like the energy is not, in general, guaranteed with conventional MOR  techniques. The violation of such invariants often results in a qualitatively wrong or unstable reduced system, even when the high-fidelity system is stable. Preservation of the conservative properties of the ROMs  by the reduced system results in  physically meaningful reduced systems for fluid dynamics problems, such as incompressible and compressible  Euler equation \cite{Afkham20}.
The stability of reduced models over long-time integration has been investigated in the context of Lagrangian systems \cite{Carlberg15}, and for port-Hamiltonian systems \cite{Chaturantabu16}.
 For linear and nonlinear Hamiltonian systems, like the linear wave equation, Sine-Gordon equation, nonlinear Schr\"odinger equation,  symplectic model reduction techniques with symplectic Galerkin projection are  constructed \cite{Hesthaven16,Peng16}  that capture the symplectic structure of Hamiltonian systems to ensure long term stability of the reduced model.

In this work, we study  structure preserving  reduced order modeling for RTSWE that exploits skew-symmetry of the centered discretization schemes to recover conservation of the energy at the level of the reduced system. Discretization of the RTSWE in space by centered finite differences leads to a skew-gradient system of ordinary differential equations (ODEs).
The skew-symmetric structure of the skew-gradient ODEs of   Hamiltonian systems  such as  Korteweg-de Vries (KdV) equation with the constant Poisson structure, i.e., skew-symmetric matrix,   are  preserved  by the ROMs using the modified POD-Galerkin
projection  \cite{Gong17,Hesthaven20,Miyatake19}.  But the skew-symmetric Poisson matrix of the RTSWE in the semi-discrete form  is state-dependent. Therefore preservation of the skew-gradient structure in the reduced form  for PDEs like the RTSWE with the state-dependent skew-symmetric Poisson structure is more challenging than with the constant Poisson structure.
  In general, there exists no numerical integration method that preserves both
the symplectic/Poisson  structure and the energy of a Hamiltonian system \cite{Hairer16}. Since the reduced-order RTSWE keeps
the skew-gradient structure and thus
have the energy-preservation law, energy-preserving integrators
can be easily applied.
The semi-discrete skew-gradient RTSWE  is solved in time with the energy preserving average vector field (AVF)  method \cite{Cohen11} which is conjugate to a Poisson integrator.
The AVF method  is  applied for  reduced order modeling of Hamiltonian systems like  the Korteweg-de Vries equation \cite{Gong17,Hesthaven20,Miyatake19} and nonlinear Schr\"odinger equation \cite{Karasozen18}. To accelerate the computation of the nonlinear terms in the reduced form, usually hyper-reduction techniques such as empirical interpolation method (EIM) and discrete empirical interpolation method (DEIM) \cite{Barrault04,chaturantabut10nmr} are used.
But the straightforward application of a structure-preserving hyper-reduction technique does not allow separation of online and offline computation of the nonlinear terms of the non-canonical Hamiltonian systems like the RTSWE, and consequently the online computational cost is not reduced. Approximation of the Poisson matrix and the gradient of the Hamiltonian of the RTSWE  with the DEIM, results in a skew-gradient reduced system with linear and quadratic terms only.
Utilizing  tensor techniques  \cite{Benner18,Benner21,Karasozen21},   the computation of the ROMs is  further accelerated while preserving the skew-symmetric structure of the RTSWE.
The cubic Hamiltonian (energy) of the RTSWE and the linear and quadratic Casimirs, i.e., the mass, the buoyancy, and the total vorticity, are preserved by  the fully discrete ROM setting applying  the POD and DEIM.
Numerical simulations for the double vortex test case from \cite{Eldred19}  confirm the structure preserving features of the ROMs,  i.e., preservation of  the Hamiltonian (energy),  mass,  buoyancy, and total potential vorticity   of the high-fidelity  FOM in long term integration, that results in the construction of a physically meaningful reduced system. Numerical results show that ROMs with the POD and DEIM  provide accurate and stable approximate solutions of the RTSWE while exhibiting several orders of magnitude computational speedup over the FOM.

The paper is organized as follows. In Section~\ref{sec:swe}, the RTSWE  is  briefly described. The Hamiltonian structure preserving FOM in space and time is introduced  in Section~\ref{sec:fom}. The ROMs with the POD and the DEIM are constructed  in Section~\ref{sec:rom}. In Section~\ref{sec:num}, numerical results  are presented. The paper ends with some conclusions in Section~\ref{sec:conc}.

%%%%%%%%%%%%%%%%%%%%%%%%%%%%%%%%%%%%%%%%%%%%%%%%%%%%%%%%%%%%%%%%%%%%%%%%%%%%%%%%%%%%%%%
%%%%%%%%%%%%%%%%%%%%%%%%%%%%%%%%%%%%%%%%%%%%%%%%%%%%%%%%%%%%%%%%%%%%%%%%%%%%%%%%%%%%%%%
%%%%%%%%%%%%%%%%%%%%%%%%%%%%%%%%%%%%%%%%%%%%%%%%%%%%%%%%%%%%%%%%%%%%%%%%%%%%%%%%%%%%%%%
\section{Rotating thermal shallow water equation}
\label{sec:swe}

RSWE  is a classical modeling tool in geophysical fluid dynamics for understandig  a variety of major dynamical phenomena in the atmosphere and oceans. But it lacks of horizontal gradients of temperature and/or density. The assumption of horizontal homogeneity of temperature/density accompanies the standard derivation of the RSWE. Relaxing this assumption  does not substantially alter the derivation and leads to the RTSWE.  Structurally, while the classical SWEs are equivalent to those of isentropic gas dynamics with pressure depending only on density, the thermal shallow water  model corresponds to the dynamics of a gas with a specific equation of state, depending both on density and temperature \cite{Zeitlin18}.

The RTSWE  known also as Ripa equation \cite{Ripa95} is obtained along the same
lines as the RSWE, by vertical averaging of the primitive
equations in the Boussinesq approximation, and using the hypothesis of
columnar motion (mean-field approximation), but relaxing the
hypothesis of uniform density/temperature \cite{Gouzien17,Zeitlin18}.
When the layer depth is supposed to be small compared
with a typical horizontal length scale,  the vertical fluid acceleration in the
layer may be neglected.
RTSWE equation is given for the primitive variables as \cite{Dellar03,Eldred19}
\begin{equation} \label{tswe}
\begin{aligned}
\frac{\partial h}{\partial t } & =   - (hu)_x   - (hv)_y\\
\frac{\partial u}{\partial t } & =  hqv -\left (\frac{u^2+v^2}{2}\right)_x- \frac{h}{2}s_x - s( h + b)_x   \\
\frac{\partial v}{\partial t } & =   -hqu -\left (\frac{u^2+v^2}{2}\right)_y - \frac{h}{2}s_y - s(h + b  )_y\\
\frac{\partial s}{\partial t } & =   - us_x   - vs_y
\end{aligned}
\end{equation}
where $u(x,y,t)$ and $v(x,y,t)$ are the relative velocities, $h(x,y,t)$ is the fluid height, $b(x,y)$ topographic height, $\rho(x,y,t)$ fluid density, $g$ gravity constant, $s = g\rho/\bar{\rho}$ the buoyancy,
$q = (v_x - u_y + f)/h$ is the potential vorticity with the constant Coriolis force $f$.
 The RTSWE equation reduces to the RSWE equation for constant buoyancy  $s$.
The RTSWE equation \eqref{tswe} is considered on a time interval $[0,T]$ for a final time $T>0$, and on a two-dimensional space domain $\Omega \in \mathbb{R}^2$
with periodic boundary conditions. The initial conditions are
$$
u({\mathbf x},0) = u_0({\mathbf x}), \; v({\mathbf x},0) = v_0({\mathbf x}), \; h({\mathbf x},0) = h_0({\mathbf x}), \;
s({\mathbf x},0) = s_0({\mathbf x}),
$$
where we set ${\mathbf x}=(x,y)^T$.
The RTSWE equation has similar  non-canonical Hamiltonian/Poisson structure as the RSWE equation \cite{Dellar03,Dellar13,Eldred19}
\begin{equation} \label{tsweh}
\dfrac{\partial z}{\partial t} = -{\mathcal J}(z)\dfrac{\delta \mathcal{H}}{\delta z}
= -
\begin{pmatrix}
0& \partial_x & \partial_y&0 \\
\partial_x& 0 & -q & -h^{-1}s_x \\
\partial_y& q  & 0 &  -h^{-1}s_y \\
0 &  h^{-1}s_x  &  h^{-1}s_y    & 0
\end{pmatrix}
\begin{pmatrix}
\frac{u^2 + v^2}{2} + sh + sb \\
hu \\
hv \\
\frac{h^2}{2} + hb
\end{pmatrix},
\end{equation}
with the Hamiltonian
\begin{equation} \label{hams}
\mathcal{H}(z)= \int_{\Omega}  \left (\frac{h^2s}{2} + hsb + h \frac{u^2 + v^2}{2} \right ) d\Omega.
\end{equation}
Substituting the variational derivatives of the  Hamiltonian \eqref{hams}
$$
\frac{\partial\mathcal{H}}{\partial {\bm u} } = h {\bm u}, \;\frac{\partial\mathcal{H}}{\partial  s } = \frac{1}{2}h^2 + hb,\;
\;\frac{\partial\mathcal{H}}{\partial  h } = \frac{1}{2}{\bm u}\cdot {\bm u} + sh + sb,\;
$$
into \eqref{tsweh} leads to the RTSWE equation \eqref{tswe} with ${\bm u}  = (u,v)^T$.\\

The non-canonical Hamiltonian form of the RTSWE equation \eqref{tsweh} is determined by the skew-symmetric Poisson bracket of two functionals $\mathcal{A}$ and $\mathcal{B}$ \cite{Salmon04} as
\begin{equation} \label{bracket}
\{ \mathcal{A},\mathcal{B}  \} = \iint \left (q\frac{\delta((\mathcal{A},\mathcal{B})}{\delta(u,v)} - \frac{\delta \mathcal{A}}{\delta {\mathbf \upsilon}}\cdot \nabla \frac{\delta \mathcal{B}}{\delta h}
+ \frac{\delta \mathcal{B}}{\delta {\mathbf \upsilon}}\cdot \nabla \frac{\delta \mathcal{A}}{\delta h} \right ) d {\mathbf x},
\end{equation}
where $\nabla =(\partial x, \partial y)^T$, and $\delta \mathcal{A}/\delta {\mathbf \upsilon}$ is the functional derivative of $\mathcal{A}$ with respect to ${\mathbf \upsilon}$. The functional Jacobian is given by
\begin{equation*}
\frac{\delta((\mathcal{A},\mathcal{B})}{\delta(u,v)} = \frac{\delta \mathcal{A}}{\delta u}\frac{\delta \mathcal{B}}{\delta v } - \frac{\delta \mathcal{B}}{\delta u}\frac{\delta \mathcal{A}}{\delta v}.
\end{equation*}
The Poisson bracket  \eqref{bracket} is related to the skew-symmetric Poisson matrix ${\mathcal J}$ as $\{\mathcal{A},\mathcal{B}\}=\{\mathcal{A},{\mathcal J}\mathcal{B}\}$.
Although the matrix  ${\mathcal J}$ in \eqref{tsweh} is not skew-symmetric, the skew-symmetry of the Poisson bracket appears after integrations by parts, and the Poisson bracket satisfies the Jacobi identity
\begin{equation*}
\{ \mathcal{A},\{\mathcal{B},\mathcal{D}\}\} + \{ \mathcal{B},\{\mathcal{D},\mathcal{A}\}\} + \{ \mathcal{A},\{\mathcal{B},\mathcal{D}\}\} = 0,
\end{equation*}
for any three functionals $\mathcal{A}$, $\mathcal{B}$ and $\mathcal{D}$.
Conservation of the Hamiltonian \eqref{hams}
follows from the antisymmetry of the Poisson bracket \eqref{bracket}
\begin{equation*}
\frac{d{\mathcal H}}{dt}= \{{\mathcal H}, {\mathcal H} \} = 0.
\end{equation*}
Other conserved quantities are the Casimirs which are additional constants of motion, and commute with any functional $\mathcal{A}$, i.e., the Poisson bracket vanishes
$$
\{\mathcal{A},\mathcal{C}\} = 0, \quad \forall \mathcal{A}({\mathbf z}) \quad \text{or} \quad \mathcal{J}^{ij}\frac{\partial \mathcal{C}}{\partial {\mathbf z}^j } =0.
$$
The Casimirs of the RTSWE equation are the  mass, the total potential vorticity, and  the  buoyancy
\begin{equation} \label{cas}
{\mathcal M} =\int h \;d\Omega, \quad {\mathcal Q }=\int hq \;d\Omega, \quad {\mathcal B} =\int hs \; d\Omega.
\end{equation}
Unlike the RSWE equation, the potential enstropy is not conserved.

%%%%%%%%%%%%%%%%%%%%%%%%%%%%%%%%%%%%%%%%%%%%%%%%%%%%%%%%%%%%%%%%%%%%%%%%%%%%%%%%%%%%%%%
%%%%%%%%%%%%%%%%%%%%%%%%%%%%%%%%%%%%%%%%%%%%%%%%%%%%%%%%%%%%%%%%%%%%%%%%%%%%%%%%%%%%%%%
%%%%%%%%%%%%%%%%%%%%%%%%%%%%%%%%%%%%%%%%%%%%%%%%%%%%%%%%%%%%%%%%%%%%%%%%%%%%%%%%%%%%%%%
\section{Full order discretization}
\label{sec:fom}

The RTSWE \eqref{tswe} is discretized by finite differences on a uniform grid in the rectangular spatial domain $\Omega =(a,b)\times(c,d)$ with the  nodes ${\bm x}_{ij} = (x_i,y_j )^T$, where $x_i=a+(i-1)\Delta x$, $y_j=c+(j-1)\Delta y$, $i,j=1,\ldots, n+1$, and spatial mesh sizes are $\Delta x=(b-a)/n$ and $\Delta y=(d-c)/n$. Using the from bottom to top/from left to right ordering of the grid nodes, the vectors of
semi-discrete state variables are given as
\begin{equation}\label{solvec}
\begin{aligned}
\bm{h}(t) &= (h_{1,1}(t),\ldots , h_{1,n}(t),h_{2,1}(t),\ldots , h_{2,n}(t), \ldots , h_{n,n}(t))^T, \\
\bm{u}(t) &= (u_{1,1}(t),\ldots , u_{1,n}(t),u_{2,1}(t),\ldots , u_{2,n}(t), \ldots , u_{n,n}(t))^T,\\
{\bm v}(t) &= (v_{1,1}(t),\ldots , v_{1,n}(t),v_{2,1}(t),\ldots , v_{2,n}(t), \ldots , v_{n,n}(t))^T,\\
{\bm s}(t) &= (s_{1,1}(t),\ldots , s_{1,n}(t),s_{2,1}(t),\ldots , s_{2,n}(t), \ldots , s_{n,n}(t))^T,
\end{aligned}
\end{equation}
where we set $w_{i,j}(t)=w(x_i,y_j,t)$ for $w=h,u,v,s$.
We note that the degree of freedom is given by $N=n^2$ because of the periodic boundary conditions, i.e., the most right grid nodes ${\bm x}_{n+1,j}$ and the most top grid nodes ${\bm x}_{i,n+1}$ are not included. The solution vector is defined by ${\bm z}(t)=(\bm{h}(t),\bm{u}(t),{\bm v}(t), {\bm s}(t)): [0,T]\mapsto \mathbb{R}^{4N}$.
Throughout the paper, we do not explicitly represent the time dependency of the semi-discrete solutions for simplicity, and we write $\bm{u}$, ${\bm v}$, $\bm{h}$, ${\bm s}$, and  ${\bm z}$.

The first order derivatives in  space are  approximated by  one dimensional centered finite differences in $x$ and $y$ directions, respectively,  and they are  extended to two dimensions utilizing the Kronecker product.
Let $\widetilde{D}_n$  be the matrix  containing the centered finite differences  under periodic boundary conditions
$$
 \widetilde{D}_n=
\begin{pmatrix}
 0& 1&  & &-1 \\
-1& 0&1 & &   \\
  & \ddots  & \ddots  &\ddots  &   \\
  &  &-1&0 &1 \\
 1&  &  & 1&0
\end{pmatrix} \in \mathbb{R}^{n\times n}.
$$
Then, on the two dimensional mesh, the centered finite difference matrices corresponding to the first order partial derivatives $\partial_x$ and $\partial_y$ are given by
$$D_x=\frac{1}{2\Delta x}\widetilde{D}_{n}\otimes I_{n}\in\mathbb{R}^{N\times N} \; ,  \quad D_y=\frac{1}{2\Delta y}I_{n}\otimes \widetilde{D}_{n}\in\mathbb{R}^{N\times N},
$$
where $\otimes$ denotes the Kronecker product, and
 $I_n$ is the identity matrices of the  size $n\times n$.

The semi-discretization of the RTSWE \eqref{tsweh} leads to a $4N$-dimensional system of Hamiltonian ODEs in skew-gradient form
\begin{equation}\label{tswehs}
\dfrac{d {\bm z}}{d t} = -J({\bm z})\nabla_{{\bm z}} H({\bm z}),
\end{equation}
with the skew-symmetric Poisson matrix and the discrete gradient Hamiltonian given by
\begin{align*}
J({\bm z}) &=
\begin{pmatrix}
0&  D_x & D_y &0 \\
D_x & 0 & -\bm{q}^d & -(\bm{h}^{-1}\circ (D_x\bm{s}))^d \\
D_y & \bm{q}^d & 0 &  -(\bm{h}^{-1}\circ (D_y\bm{s}))^d\\
0 & (\bm{h}^{-1}\circ (D_x\bm{s}))^d &  (\bm{h}^{-1}\circ (D_y\bm{s}))^d & 0
\end{pmatrix} ,\\[0.3cm]
\nabla_{{\bm z}} H({\bm z}) &=
\begin{pmatrix}
\frac{\bm{u}^2 + {\bm v}^2}{2} + {\bm s}  \circ \bm{h} + b{\bm s} \\
\bm{h} \circ  \bm{u}\\
\bm{h} \circ  {\bm v}\\
\frac{1}{2}\bm{h}^2 + b\bm{h}
\end{pmatrix},
\end{align*}
where $ \circ $ denotes element-wise (or Hadamard) product, and the square operations are also held element-wise. The matrix  ${\bm q}^d\in \mathbb{R}^{N\times N}$ is the diagonal matrix with the diagonal elements $({\bm q}^d)_{ii}={\bm q}_i$ where ${\bm q}$ is the semi-discrete vector of the potential vorticity $q(x,y,t)$, $i=1,\ldots,N$. The diagonal matrices with the diagonal elements from the vectors $\bm{h}^{-1}\circ (D_x\bm{s})$ and $\bm{h}^{-1}\circ (D_y\bm{s})$ are defined similarly.

We use the energy  preserving AVF method \cite{Cohen11} for solving semi-discretized RTSWE \eqref{tswehs} in time. The AVF method is used as the time integrator for the RSWE \cite{Cotter18,Wimmer20a} and for the linearized RTSWE \cite{Eldred19} in non-canonical Hamiltonian form. The AVF method  preserves exactly the Hamiltonian (energy) and the quadratic Casimirs of the non-canonical Hamiltonian systems in skew-gradient (Poisson) form such as the semi-discretized RTSWE \eqref{tswehs}. It is invariant with respect to linear transformations, symmetric and of order 2. The Hamiltonian  (energy) and the symplectic or Poisson structure can not be preserved simultaeously   by any time integrator  \cite{Hairer16}. But, the energy preserving AVF method is conjugate to a Poisson integrator in the sense that all Casimir functions of the RTSWE are nearly conserved without drift as will be shown  in  Section \ref{sec:num}. Linear conserved quantities like the total mass is preserved by all time integrators including the AVF method.

The time  interval $[0,T]$ is partitioned into $K$ uniform meshes with the step size $\Delta t=T/K$ as $0=t_0<t_1<\ldots <t_K=T$, and $t_k=k\Delta t$, $k=0,1,\ldots ,K$.  The fully discrete solution vector at time $t_k$ is denoted as $\bm{z}^k=\bm{z}(t_k)$.  Similar setting is used for the other state variables. Integration of the semi-discrete RTSWE \eqref{tswehs} in time  by the AVF integrator leads to the fully discrete problem: for $k=0,1,\ldots , K-1$, given ${\bm z}^{k}$ find ${\bm z}^{k+1}$ satisfying
\begin{align}\label{tswe_avf}
 {\bm z}^{k+1}&={\bm z}^{k}+ \Delta t J\left (\dfrac{{\bm z}^{k+1}+{\bm z}^{k}}{2}\right )\int_{0}^{1} \nabla_{\bm z} H(\xi({\bm z}^{k+1}-{\bm z}^{k})+{\bm z}^{k})d \xi.
\end{align}

Practical implementation of the AVF method requires the evaluation of the integral on the right-hand side of \eqref{tswe_avf}. Since the Hamiltonian ${\mathcal H}$ and Casimirs are polynomial, they can be exactly integrated with the Gaussian quadrature rule of the appropriate degree. Since
the AVF method is an implicit time integrator, at each time step, a non-linear system
of equations has to be solved by the Newton's method.

Under periodic boundary conditions, the full discrete forms of the energy $H(\bm{z})$, the
 mass $M(\bm{z})$, the total potential vorticity $Q(\bm{z})$, and the buoyancy $B(\bm{z})$ are given at the time instance $t_k$ as
\begin{equation}\label{dener}
\begin{aligned}
H^k(\bm{z}) &= \sum_{i=1}^n\sum_{j=1}^n \left (\frac{1}{2}(h_{i,j}^k)^2 s_{i,j}^k  + h_{i,j}^k s_{i,j}^k b_{i,j}^k + h_{i,j}^k\frac{(u_{i,j}^k)^2 + (v_{i,j}^k)^2}{2}   \right ) \Delta x\Delta y,
\\
M^k(\bm{z}) &= \sum_{i=1}^n\sum_{j=1}^n  h_{i,j}^k \Delta x\Delta y, \\
Q^k(\bm{z}) & =  \sum_{i=1}^n\sum_{j=1}^n \left ( \frac{v_{i+1,j}^k-v_{i-1,j}^k}{2\Delta x} -  \frac{u_{i,j+1}^k-u_{i,j-1}^k}{2\Delta y} + f \right ) \Delta x\Delta y, \\
B^k(\bm{z}) &= \sum_{i=1}^n\sum_{j=1}^n h_{i,j}^k s_{i,j}^k \Delta x\Delta y,
\end{aligned}
\end{equation}
where $b_{i,j}^k\approx b(x_i,y_j,t_k)$ and for $w=h,u,v,s$, $w_{i,j}^k\approx h(x_i,y_j,t_k)$ are full discrete approximations at $(x_i,y_j,t_k)$ satisfying  $w_{i+n,j}^k=w_{i,j}^k$  and $w_{i,j+n}^k=w_{i,j}^k$  because of the periodic boundary conditions, $i,j=1,\ldots ,n$.

We remark that recently the RTSWE was solved with the well-balanced finite volume method as hyperbolic system in \cite{Kurganov20,Kurganov20a}.

%%%%%%%%%%%%%%%%%%%%%%%%%%%%%%%%%%%%%%%%%%%%%%%%%%%%%%%%%%%%%%%%%%%%%%%%%%%%%%%%%%%%%%%
%%%%%%%%%%%%%%%%%%%%%%%%%%%%%%%%%%%%%%%%%%%%%%%%%%%%%%%%%%%%%%%%%%%%%%%%%%%%%%%%%%%%%%%
%%%%%%%%%%%%%%%%%%%%%%%%%%%%%%%%%%%%%%%%%%%%%%%%%%%%%%%%%%%%%%%%%%%%%%%%%%%%%%%%%%%%%%%
\section{Reduced order model}
\label{sec:rom}

In this section, we construct ROMs that preserve the skew-gradient structure of the semi-discrete RTSWE \eqref{tswehs}, and consequently the discrete conserved quantities in \eqref{dener}.
Because the RTSWE is a non-canonical Hamiltonian PDE with a state-dependent Poisson structure, a straightforward application of the POD will not preserve the skew-gradient structure of the RTSWE \eqref{tswehs} in the reduced form.
The AVF method was used as an energy  preserving time integrator for reduced order modelling of  Hamiltonian systems with constant skew-symmetric matrices like the Korteweg de Vires equation  \cite{Gong17,Miyatake19} and nonlinear Schr\"odinger equation
(NLSE) \cite{Karasozen18}.  The approach in \cite{Gong17} can  also be  applied to skew-gradient systems with state-dependent skew-symmetric structure as the RTSWE equation \eqref{tswehs}.
We show that the state-dependent skew-symmetric matrix in \eqref{tswehs} can be evaluated efficiently in the online stage independent of the full dimension $N$. The full and reduced models are computed separately approximating the nonlinear terms by DEIM.

The POD basis vectors are obtained usually by stacking all the state variables in one vector and  a  common reduced subspace is computed by taking the singular value decomposition (SVD) of the snapshot data. Because the governing PDEs like the RTSWE are coupled, the resulting ROMs do not preserve the coupling structure of the FOM  and  produce unstable ROMs  \cite{Kramer20,Reiss07}. In order to maintain the coupling structure in the ROMs, the POD basis vectors are computed separately for each the state vector $\bm{h}$, $\bm{u}$, $\bm{v}$ and $\bm{s}$.

The POD basis is computed through the mean subtracted snapshot matrices ${\mathsf S}_u$, ${\mathsf S}_v$, ${\mathsf S}_h$  and ${\mathsf S}_s$, constructed by the solutions of the full discrete high fidelity model \eqref{tswe_avf}
\begin{align*}
{\mathsf S}_u &= \left({\bm u}^1 - \overline{\bm u},  \cdots, {\bm u}^{K} - \overline{\bm u} \right) \in\mathbb{R}^{N\times K}, \\
{\mathsf S}_v &= \left({\bm v}^1 - \overline{\bm v}, \cdots, {\bm v}^{K} - \overline{\bm v} \right)   \in\mathbb{R}^{N\times K},  \\
{\mathsf S}_h &= \left({\bm h}^1 - \overline{\bm h}, \cdots, {\bm h}^{K}  - \overline{\bm h}\right) \in\mathbb{R}^{N\times K},\\
{\mathsf S}_s &= \left({\bm s}^1 - \overline{\bm s}, \cdots, {\bm s}^{K}  - \overline{\bm s}\right) \in\mathbb{R}^{N\times K},
\end{align*}
where $\overline{\bm u}$, $\overline{\bm v}$, $\overline{\mathbf{h}}$, $\overline{\mathbf{s}}\in\mathbb{R}^{N}$  denote the time averaged means of the solutions
$$
\overline{\bm u} = \frac{1}{K}\sum_{k=0}^{K} {\bm u}^k\; , \; \overline{\bm v} = \frac{1}{K}\sum_{k=0}^{K} {\bm v}^k\; , \; \overline{\bm h} = \frac{1}{K}\sum_{k=0}^{K} { \bm h}^k, \; \overline{\bm s} = \frac{1}{K}\sum_{k=0}^{K} { \bm s}^k
$$
The mean-subtracted ROMs is used frequently in fluid dynamics to stabilize the reduced system, and it guarantees that ROM solutions would
satisfy the same boundary conditions as for the FOM \cite{Navon20,Bistrian15}.

The POD modes are computed by applying SVD to the snapshot matrices as
\begin{equation*}
{\mathsf S}_u= W_u \Sigma_u U_u^T\, \quad {\mathsf S}_v= W_v \Sigma_v U_v^T, \quad {\mathsf S}_h= W_h \Sigma_h U_h^T, \quad {\mathsf S}_s= W_s \Sigma_s U_s^T,
\end{equation*}
where for $i=u,v,h,s$, $W_i \in \mathbb{R}^{ N\times K}$ and
$U_i\in \mathbb{R}^{ K\times K}$ are orthonormal matrices,  and $\Sigma_i\in \mathbb{R}^{ K\times K}$ is the diagonal matrix with its diagonal entries are the singular values $\sigma_{i,1} \ge \sigma_{i,2} \ge \cdots \ge \sigma_{i,K}\geq 0$.
The matrices $V_{i,r_i}\in\mathbb{R}^{N\times r_i}$ of rank $r_i$ POD modes consists of the first $r_i$ left singular vectors from $W_i$ corresponding to the $r_i$ largest singular values, which satisfies the following least squares error
$$
\min_{V_{i,r_i}\in \mathbb{R}^{ N\times r_i}}||{\mathsf S}_i-V_{i,r_i}V_{i,r_i}^T{\mathsf S}_i ||_F^2 = \sum_{j=r_i+1}^{K} \sigma_{i,j}^2\; , \quad i=u,v,h,s,
$$
where $\|\cdot\|_F$ is the Frobenius norm. In the sequel, we omit the subscript $r_i$ from the rank $r_i$ POD matrices $V_{i,r_i}$, and we write $V_{i}$
for the ease of notation. Moreover, we have the reduced approximations
\begin{equation}\label{relz}
{\bm u} \approx  \widehat{{\bm u}}=\overline{\bm u} + V_{ u}{\bm u}_r, \quad {\bm v} \approx  \widehat{{\bm v}}=\overline{\bm v} + V_{v}{\bm v}_r,
\quad {\bm  h} \approx \widehat{{\bm h}}=\overline{\bm h} + V_{h}  {\bm h}_r,
\quad {\bm  s} \approx \widehat{{\bm s}}=\overline{\bm s} + V_{s}  {\bm s}_r,
\end{equation}
where the  vectors ${\bm u}_r$, ${\bm v}_r,{\bm h}_r $ and ${\bm s}_r$ are the ROM solutions which are the coefficient vectors for the reduced approximations $\widehat{{\bm u}}$, $\widehat{{\bm v}}$, $\widehat{{\bm h}}$ and $\widehat{{\bm s}}$. The ROM is  obtained by Galerkin projection onto the reduced space
\begin{align} \label{galpod1}
\frac{d}{dt}{\bm z}_r = V_{z}^T J(\widehat{{\bm z}})\nabla_{\bm z} H(\widehat{{\bm z}}),
\end{align}
where the ROM solution is ${\bm z}_r=({\bm h}_r,{\bm u}_r,{\bm v}_r, {\bm s}_r): [0,T]\mapsto \mathbb{R}^{r_h+r_u+r_v+r_s}$ and the related reduced approximation is $\widehat{{\bm z}}=(\widehat{{\bm h}},\widehat{{\bm u}},\widehat{{\bm v}}, \widehat{{\bm s}}) $.
The block diagonal matrix $V_{z}$ contains the matrix of POD modes for each state variable
\begin{equation*}
V_{z}=
\begin{pmatrix}
V_{h} &  & &\\
& V_{u}& & \\
& & V_{v}  & \\
& & & V_{s}   \\
\end{pmatrix}\in\mathbb{R}^{4N\times (r_h+r_u+r_v+r_s)}.
\end{equation*}
Although the formulae  in this section are generally applicable for any
number of modes $r_h,\; r_u,\; r_v$ and $r_s$, we use equal number of POD modes in the present study  $r := r_h=r_u = r_v = r_s$.

The skew-gradient structure of the FOM \eqref{tswehs} is  not preserved  by the ROM \eqref{galpod1}.
A reduced skew-gradient system is obtained formally by inserting $V_{z}V_{z}^T$ between $J(\widehat{{\bm z}})$ and $\nabla_{\bm z} H(\widehat{{\bm z}})$  \cite{Gong17,Miyatake19}, leading to the ROM
\begin{align} \label{galpod}
\frac{d}{dt}{\bm z}_r = J_r(\widehat{{\bm z}})\nabla_{\mathbf{z}_r} H(\widehat{{\bm z}}),
\end{align}
where $ J_r(\widehat{{\bm z}})= V_{z}^T J( \widehat{{\bm z}})V_{z}$ and $\nabla_{\mathbf{z}_r} H(\widehat{{\bm z}})= V_{z}^T\nabla_{\mathbf{z}} H(\widehat{{\bm z}})$.
The reduced order RTSWE  \eqref{galpod} is also solved by the AVF method.

Using the above approach, in \cite{Gong17,Miyatake19} non-canonical Hamiltonian systems with constant skew-symmetric matrices, like the Korteweg de Vries equations are solved with the POD and DEIM. As time integrator the AVF method and the mid-point rule is used. Because the skew-symmetric matrix of the semi-discretized RTSWE
\eqref{tswehs} depends on the stare variables  nonlinearly, the construction of the ROMs while preserving the skew-gradient structure in \eqref{galpod} is more challenging. In the following we describe structure-preserving POD and DEIM for the semi-discretized RTSWE \eqref{tswehs} in the reduced form \eqref{galpod}

The skew-symmetric matrix $J_r(\widehat{{\bm z}})$ of the reduced order RTSWE  \eqref{galpod} is given explicitly as
\begin{equation*}
J_r(\widehat{{\bm z}}) =
\begin{pmatrix}
0 & -V_{h}^TD_xV_{u} & - V_{h}^TD_yV_{v}  &  0\\
-V_{u}^TD_xV_{h} &  0  & V_{u}^TF_1(\bm{\widehat{z}})^{d}V_{v}    & V_{u}^TF_2(\bm{\widehat{z}})^dV_{s}\\
- V_{v}^TD_yV_{h} & -V_{v}^TF_1(\bm{\widehat{z}})^{d}V_{u}  &  0  & V_{v}^TF_3(\bm{\widehat{z}})^dV_{s}\\
0  &  -V_{s}^TF_2(\bm{\widehat{z}})^dV_{u}  & -V_{s}^TF_3(\bm{\widehat{z}})^dV_{v} & 0,
\end{pmatrix}
\end{equation*}
where
$$
F_1(\bm{\widehat{z}})=\widehat{{\bm q}}, \quad F_2(\bm{\widehat{z}})=\widehat{{\bm h}}^{-1}\circ (D_x\widehat{{\bm s}}),\quad F_3(\bm{\widehat{z}})=
\widehat{{\bm h}}^{-1}\circ (D_y\widehat{{\bm s}})
$$
and the superscript $ d $ refers to the  diagonal matrix of the vector as in \eqref{tswehs}.
The matrices $F_j(\bm{\widehat{z}})$, $j=1,2,3$,  are not constant and they should be computed in the online stage depending on the full order system, whereas other matrices  are constant and they can be  precomputed in the offline stage.

The semi-discrete structure-preserving ROM \eqref{galpod} of the \eqref{tswehs} can be written explicitly as
\begin{equation}\label{discpod}
\begin{aligned}
\frac{d}{dt}{\bm h}_r&= - (V_{h}^TD_xV_{u}V_{u}^T{F_5} + V_{h}^TD_yV_{v}V_{v}^T{F_6})\\
\frac{d}{dt}{\bm u}_r&= - (V_{u}^TD_xV_{h}V_{h}^T{F_4} - V_{u}^T({F}_1\circ V_{v}V_{v}^T{F_6})-V_{u}^T({F}_2\circ V_{s}V_{s}^T{F_7}))\\
\frac{d}{dt}{\bm v}_r&= - (V_{v}^TD_yV_{h}V_{h}^T{F_4} + V_{v}^T({F}_1\circ V_{u}V_{u}^T{F_5}) - V_{v}^T({F}_3\circ V_{s}V_{s}^T{F_7}))\\
\frac{d}{dt}{\bm s}_r&= - (V_{s}^T({F}_2\circ V_{u}V_{u}^T{F_5}) + V_{s}^T({F}_3\circ V_{v}V_{v}^T{F_6})),
\end{aligned}
\end{equation}
where
$$F_4(\bm{\widehat{z}})=\frac{\bm{\widehat{u}}^2 + \bm{\widehat{v}}^2}{2} + \bm{\widehat{s}}  \circ \bm{\widehat{h}} + b\bm{ \widehat{s}}, \;F_5(\bm{\widehat{z}}) = \bm{\widehat{h}} \circ  \bm{\widehat{u}}, \;
F_6(\bm{\widehat{z}}) =   \bm{\widehat{h}} \circ  \bm{\widehat{v}}, \;
    F_7(\bm{\widehat{z}}) =  	\frac{\bm{\widehat{h}}^2}{2} + b\bm{\widehat{h}}.
$$

Due to the nonlinear terms, the computation of the reduced system \eqref{discpod} still scales with the dimension $N$ of the FOM. This can be handled   by applying a hyper-reduction technique  such as the DEIM \cite{chaturantabut10nmr}. To preserve the structure in \eqref{galpod}, we approximate the nonlinear terms with DEIM for skew-symmetric reduced Poisson matrix $ J_r(\widehat{{\bm z}}) $ and discrete reduced gradient of the Hamiltonian $ \nabla_{\mathbf{z}_r} H(\widehat{{\bm z}}) $ separately.

The DEIM procedure, originally introduced in \cite{chaturantabut10nmr}, is utilized to approximate the nonlinear vectors $F_j(\widehat{{\bm z}}(t))$ in \eqref{discpod}
by interpolating them onto an empirical basis, that is,
$$
F_j(\widehat{{\bm z}}) \approx \Phi^{(j)} \mathsf{c}^{(j)}(t), \quad j=1,\ldots ,7,
$$
where $\Phi^{(j)} = [\phi_{1}^{(j)}, \dots, \phi_{p_j}^{(j)}] \in \mathbb{R}^{N \times p_j}$ is a low dimensional basis matrix
and $\mathsf{c}^{(j)}(t) : [0,T] \mapsto \mathbb{R}^{p_j}$ is the vector of time-dependent coefficients to be determined.
Let
$\mathsf{P}^{(j)} = [e_{\rho_{1}}^{(j)}, \dots, e_{\rho_{p}}^{(j)}] \in \mathbb{R}^{N \times p_j}$ be a subset of columns of the identity
matrix, named as the "selection matrix".
If $(\mathsf{P}^{(j)})^{\top}\Phi^{(j)}$ is invertible,
in \cite{chaturantabut10nmr} the coefficient vector $\mathsf{c}^{(j)}(t)$ is uniquely determined by solving the linear system
$(\mathsf{P}^{(j)})^{\top}\Phi^{(j)} \mathsf{c}^{(j)}(t) = (\mathsf{P}^{(j)})^{\top}F_j(\widehat{{\bm z}}(t))$,
so that the nonlinear terms in the reduced model \eqref{discpod} is then approximated by
\begin{equation} \label{fkapprox}
F_j(\widehat{{\bm z}}(t)) \approx \Phi^{(j)} \mathsf{c}^{(j)}(t) = \Phi^{(j)}((\mathsf{P}^{(j)})^{\top}\Phi^{(j)})^{-1}(\mathsf{P}^{(j)})^{\top}F_j(\widehat{{\bm z}}(t)).
\end{equation}

The accuracy of DEIM depends mainly  on the choice of the  basis, and not much by the choice of $\mathsf{P}^{(j)}$.
In most applications,  the interpolation basis $\{\phi_1^{(j)},\dots,\phi_p^{(j)}\}$ is selected as the POD basis of the set of snapshots matrices
of the nonlinear vectors given by
\begin{equation}
\mathsf{S}_{F_j} = [ F_j^1,F_j^2,\cdots , F_j^{K}]\in\mathbb{R}^{N\times K}, \quad j=1,\ldots ,7,
\label{nonsvd}
\end{equation}
where $F_j^k=F_j(\widehat{{\bm z}}^k)$ denotes the nonlinear vector $F_j(\widehat{{\bm z}}(t))$ at time $t_k$, computed by using the solution vectors $\widehat{{\bm z}}^k=\widehat{{\bm z}}(t_k)$, $k=1,\ldots , K$.
The columns of the matrices $\Phi^{(j)} = [\phi_{1}^{(j)}, \dots, \phi_{p_j}^{(j)}]$ are determined as the first $p_j\ll N$ dominant left singular vectors in the SVD of $\mathsf{S}_{F_j}$.
The selection matrix $\mathsf{P}^{(j)}$ for DEIM is determined by a greedy algorithm based on the system residual;
see \cite[Algorithm 3.1]{chaturantabut10nmr}. Similar to the POD, we take the same  number of DEIM modes for the non linear terms, i.e., $p:=p_j, j=1,\ldots 7$.

An elegant way of approaching the sampling points  is the Q-DEIM \cite{gugercin16} which relies on QR decomposition with column pivoting. It was shown that Q-DEIM  leads to better accuracy and stability properties of
the computed selection matrix $\mathsf{P}^{(j)}$ with the pivoted
QR-factorization of $(\Phi^{(j)})^{\top}$. In the sequel, we use Q-DEIM for the calculation of the selection matrix $\mathsf{P}^{(j)}$, see Algorithm~\ref{alg:qdeim}.

\begin{algorithm}[H]
\caption{Q-DEIM \label{alg:qdeim}}
\begin{algorithmic}[1]
	\STATE\textbf{Input:} Basis matrix $\Phi\in\mathbb{R}^{N\times p}$
	\STATE\textbf{Output:} Selection matrix $\mathsf{P}$
	\STATE Perform pivoted QR factorization of $\Phi^{\top}$ so that $\Phi^{\top}\Pi = QR$
	\STATE Set $\mathsf{P} = \Pi (:,1:p)$
\end{algorithmic}
\end{algorithm}

Applying the  DEIM/Q-DEIM  \eqref{fkapprox},
the  skew-gradient structure of the reduced system of the form
\begin{align} \label{hamDEIM}
\frac{d}{dt}{\bm z}_r = \widetilde{J_r}(\widehat{{\bm z}})\nabla_{\mathbf{z}_r} \widetilde{H}(\widehat{{\bm z}}),
\end{align}

\begin{equation*}
\widetilde{J_r}(\widehat{{\bm z}}) =
\begin{pmatrix}
0 & -V_{h}^TD_xV_{u} & - V_{h}^TD_yV_{v}  &  0\\
-V_{u}^TD_xV_{h} &  0  & V_{u}^T\widetilde{F_1}(\bm{\widehat{z}})^{d}V_{v}    & V_{u}^T\widetilde{F_2}(\bm{\widehat{z}})^dV_{s}\\
- V_{v}^TD_yV_{h} & -V_{v}^T\widetilde{F_1}(\bm{\widehat{z}})^{d}V_{u}  &  0  & V_{v}^T\widetilde{F_3}(\bm{\widehat{z}})^dV_{s}\\
0  &  -V_{s}^T\widetilde{F_2}(\bm{\widehat{z}})^dV_{u}  & -V_{s}^T\widetilde{F_3}(\bm{\widehat{z}})^dV_{v} & 0
\end{pmatrix}
 \end{equation*}
is preserved, where
 $ \nabla_{\mathbf{z}_r} \widetilde{H}(\widehat{{\bm z}}) $  and $\widetilde{J_r}$ are  the DEIM reduced discrete gradient of Hamiltonian and of the  Poisson matrix, respectively and
 $$
\widetilde{F_j}=\Psi_j F_{r,j}(\widehat{{\bm z}}),\; F_{r,j}(\widehat{{\bm z}})=(\mathsf{P}^{(j)})^TF_j(\widehat{{\bm z}}),\;  \Psi_{j} =\Phi^{(j)}((\mathsf{P}^{(j)})^T\Phi^{(j)})^{-1}, \;j =1,\ldots,7.
$$

Finally the reduced system can be written as linear-quadratic ODE in the tensor form as
 \begin{equation}\label{deimtrom}
 \begin{aligned}
\frac{d}{dt}{\bm h}_r&= - \left(V_{h}^TD_xV_{u}V_{u}^T\widetilde{F_5} + V_{h}^TD_yV_{v}V_{v}^T\widetilde{F_6}\right)\\
\frac{d}{dt}{\bm u}_r&= - \left(V_{u}^TD_xV_{h}V_{h}^T\widetilde{F_4} - H_1(F_{r,1}(\widehat{{\bm z}})\otimes F_{r,6}(\widehat{{\bm z}}))-H_2(F_{r,2}(\widehat{{\bm z}})\otimes F_{r,7}(\widehat{{\bm z}}))\right)\\
\frac{d}{dt}{\bm v}_r&= - \left(V_{v}^TD_yV_{h}V_{h}^T\widetilde{F_4} + H_3(F_{r,1}(\widehat{{\bm z}})\otimes F_{r,5}(\widehat{{\bm z}})) - H_4(F_{r,3}(\widehat{{\bm z}})\otimes F_{r,7}(\widehat{{\bm z}}))\right)\\
\frac{d}{dt}{\bm s}_r&= - \left(H_5(F_{r,2}(\widehat{{\bm z}})\otimes F_{r,5}(\widehat{{\bm z}})) + H_6(F_{r,3}(\widehat{{\bm z}})\otimes F_{r,6}(\widehat{{\bm z}}))\right),
 \end{aligned}
 \end{equation}
where
\begin{equation}
\begin{aligned}
G_1&=V_{u}^TG(\Psi_{1}\otimes V_{v}V_{v}^T\Psi_{6}), \quad G_2=V_{u}^TG(\Psi_{2} \otimes V_{s}V_{s}^T\Psi_{7}),\\
G_3&=V_{v}^TG(\Psi_{1}\otimes V_{u}V_{u}^T\Psi_{5}), \quad G_4=V_{v}^TG(\Psi_{3} \otimes V_{s}V_{s}^T\Psi_{7}),\\
G_5&=V_{s}^T G(\Psi_{2}\otimes V_{u}V_{u}^T\Psi_{5}),\quad G_6=V_{s}^T G(\Psi_{3}\otimes V_{v}V_{v}^T\Psi_{6}),
\end{aligned}
\end{equation}
are the reduced matricized tensors of size $r\times p^2$, where $G\in\mathbb{R}^{N\times N^2}$ is the matricized tensor satisfying the identity $G(\bm{a}\otimes\bm{b})=\bm{a}\circ\bm{b}$ for any vectors $\bm{a},\bm{b}\in\mathbb{R}^{N}$.
All the matrices $G_j$ can be  precomputed in the offline stage, and
the reduced nonlinearities $F_{r,j}(\widehat{{\bm z}})=(\mathsf{P}^{(j)})^TF_j(\widehat{{\bm z}})$ are computed by considering just $p\ll N$ entries of the nonlinearities $F_j(\widehat{{\bm z}})$ among $N$ entries, $j=1,\ldots,7$ and $ i=h,u,v,s $. Thus, in the online phase, the cost of ROM \eqref{deimtrom} scales with $ \mathcal{O}(rp^2) $, which depends only on the reduced dimension. We remark that it is also possible to obtain a structure-preserving ROM by using only DEIM applied to the nonlinearities in the reduced Poisson matrix $ J_r $ and approximating reduced gradient of the Hamiltonian by POD. Nevertheless, in this case, the online computational cost increases rapidly  due to the  cubic nonlinear terms.

A straightforward computation of matricized tensors $ G_j $ for $ j=1,\ldots,6 $ creates big burden on the offline cost. Recently tensorial algorithms are developed by exploiting the particular structure of Kronecker product  and using the CUR matrix approximation of $ G_j $ \cite{Benner18,Benner21} to increase computational efficiency.  To reduce the computational cost of the matricized tensors $ G_j $, we follow the same procedure as in \cite{Karasozen21}. We illustrate this for  $ C=G(\Psi_{1}\otimes D) \in \mathbb{R}^{N\times p^2 }$ and $ D=V_{v}V_{v}^T\Psi_{6} \in \mathbb{R}^{N\times p}$, the remaining ones can be handled in a similar way.

Due to the large FOM dimension $N$, computational complexity of $ (\Psi_{1}\otimes D) $ is of order $ \mathcal{O}( N^2 p^2) $, which is infeasible in terms of computation and memory. The computational cost and memory can be reduced using the sparse structure of the matricized tensor $ G \in \mathbb{R}^{N\times N^2}$. The matrix $ C $ can be computed using transposed Khatri-Rao product  \cite{Benner18,Benner19} as follows
\begin{align}\label{face-splitting}
C=G(\Psi_{1}\otimes D)  =
\begin{pmatrix}
\Psi_{1}(1,:)\otimes D(1,:)\\
\vdots\\
\Psi_{1}(N,:)\otimes D(N,:)
\end{pmatrix}.
\end{align}
To increase computational efficiency in \eqref{face-splitting}, we make use of "MULTIPROD" package \cite{Leva08mmm}, which performs multiplications of the multi-dimensional arrays by applying   virtual array expansion . The Kronecker product of any two column vectors $\mathbf{a}$ and $\mathbf{b}$ satisfies
\begin{equation}\label{vec}
\begin{aligned}
\mathbf{a}^\top\otimes \mathbf{b}^\top=(\text{vec}{(\mathbf{b}\mathbf{a}^\top)})^\top,
\end{aligned}
\end{equation}
where vec$(\cdot)$ denotes the vectorization of a matrix. Then, the matrix $C$ can be constructed using \eqref{vec} as
\begin{equation}\label{id}
C(m,:)=(\text{vec}(D(m,:)^\top \Psi_{1}(m,:))^\top , \ \ m\in\{1,2,\ldots,N\}.
\end{equation}
Construction of $ C $ still requires loop operations over FOM dimension $ N $ in \eqref{id}, which can be less efficient when $ N $ is large. The multiplications \eqref{id} can be done by adding a singleton to  $\Psi_{1} \in \mathbb{R}^{ N\times 1 \times p}$ and computing MULTIPROD of $D$ and $\Psi_{1}$ in $2$nd and $3$th dimensions as
$$
\mathcal{C} =\text{MULTIPROD}(D,\Psi_{1})\in \mathbb{R}^{N\times p \times p},
$$
where the matricized tensor $C$ is obtained by reshaping the 3-dimensional array $\mathcal{C}$ into a matrix of dimension $\mathbb{R}^{N\times p^2}$.

We remark that Poisson structure or energy preserving reduced models are constructed for non-canonical Hamiltonian systems,  where the Poisson matrices depend quadratically on the states, for the Volterra lattice equation \cite{Hesthaven20} and for the modified KDV equation \cite{Miyatake19}. In both papers, tensor techniques are not exploited with the POD ad DEIM. In \cite{Miyatake19}, the Poisson matrix is approximated by DEIM but not the gradient of the Hamiltonian, which is linear in the state variable for the modified KdV equation. Our approach approximating  both the Poisson matrix and the gradient of the Hamiltonian  with the DEIM, can be applied to any non-canonical Hamiltonian system, that leads to a skew-gradient reduced system with linear and quadratic terms only. Application of the tensor techniques to the linear-quadratic reduced systems, speeds up the computation of the ROMs further.

%%%%%%%%%%%%%%%%%%%%%%%%%%%%%%%%%%%%%%%%%%%%%%%%%%%%%%%%%%%%%%%%%%%%%%%%%%%%%%%%%%%%%%%
%%%%%%%%%%%%%%%%%%%%%%%%%%%%%%%%%%%%%%%%%%%%%%%%%%%%%%%%%%%%%%%%%%%%%%%%%%%%%%%%%%%%%%%
%%%%%%%%%%%%%%%%%%%%%%%%%%%%%%%%%%%%%%%%%%%%%%%%%%%%%%%%%%%%%%%%%%%%%%%%%%%%%%%%%%%%%%%
\section{Numerical results}
\label{sec:num}

We consider  the double vortex test case from \cite{Eldred19} in the doubly periodic space domain $\Omega =[0, L]^2$ without the bottom topography $ (b =0) $. The initial conditions are given by
\begin{align*}
&h=H_0-\Delta h\left[e^{-0.5((x'_1)^2+(y'_1)^2)}+e^{-0.5((x'_2)^2+(y'_2)^2)}-\frac{4\pi \sigma_x \sigma_y}{L^2}\right],\\
&u=\frac{-g\Delta h}{f \sigma_y}\left[y_1''e^{-0.5((x'_1)^2+(y'_1)^2)}+y''_2e^{-0.5((x'_2)^2+(y'_2)^2)}\right],
\\
&v=\frac{g\Delta h}{f \sigma_x}\left[x_1''e^{-0.5((x'_1)^2+(y'_1)^2)}+x''_2e^{-0.5((x'_2)^2+(y'_2)^2)}\right],
\\
&s=g\left(1+0.05sin\left[\frac{2\pi}{L}(x-xc)\right]\right)
\end{align*}
where $ xc=0.5L $ and
\begin{align*}
	&x'_1=\frac{L}{\pi \sigma_x}\sin\left[\frac{\pi}{L}(x-xc_1)\right],\qquad x'_2=\frac{L}{\pi \sigma_x}\sin\left[\frac{\pi}{L}(x-xc_2)\right],\\
	&y'_1=\frac{L}{\pi \sigma_y}\sin\left[\frac{\pi}{L}(y-yc_1)\right],\qquad y'_2=\frac{L}{\pi \sigma_y}\sin\left[\frac{\pi}{L}(y-yc_2)\right]\\
	&x''_1=\frac{L}{2\pi \sigma_x}\sin\left[\frac{2\pi}{L}(x-xc_1)\right],\qquad x''_2=\frac{L}{2\pi \sigma_x}\sin\left[\frac{2\pi}{L}(x-xc_2)\right],\\
	&y''_1=\frac{L}{2\pi \sigma_y}\sin\left[\frac{2\pi}{L}(y-yc_1)\right],\qquad y''_2=\frac{L}{2\pi \sigma_y}\sin\left[\frac{2\pi}{L}(y-yc_2)\right]
\end{align*}
The center of the two vortices are  given by
\begin{align*}
xc_1=(0.5-ox)L,\quad xc_2=(0.5+ox)L, \quad
yc_1=(0.5-oy)L,\quad yc_2=(0.5+oy)L.
\end{align*}
The parameters are $ L=5000 $km, $  f=0.00006147\text{s}^{-1} $, $ H_0=750 $m,
$ \Delta h =75 $m, $ g=9.80616 $ms$ ^{-2} $, $ \sigma_x=\sigma_y= \frac{3}{40}L$ and $ ox =oy =0.1$.  The simulations are performed on the grid $\Omega =[0, L]^2$ with the  mesh sizes  $\Delta x = \Delta y  = 50$km, and for the number of time steps $K=250$ with the time step-size $\Delta t= 486$s, which leads to the final time $T=33$h $45$min. The size of the  snapshot matrices for each state variable is $10000\times 250$.

The POD and  DEIM  basis are truncated according to the following   relative cumulative energy criterion
\begin{equation} \label{energy_criteria}
\min_{1\leq p \leq K}\frac{\sum_{j=1}^p \sigma_{j}^2}{\sum_{j=1}^{K} \sigma_{j}^2  } > 1 - \kappa,
\end{equation}
where $\kappa $ is a user-specified tolerance.
Approximation of the nonlinear terms using DEIM can affect the accuracy of conserved quantities. Therefore, the accuracy of
the DEIM approximation must catch much more relative cumulative energy than the one for POD.
In our simulations, we set $\kappa = 10^{-3}$ and $\kappa = 10^{-5}$
to catch at least $99.9 \%$ and $99.999 \%$ of relative cumulative energy for POD and DEIM, respectively.

The accuracy of the ROMs is measured by the time averaged relative $L_2$-errors  between FOM and ROM solutions  for each state variable  ${\bm w}= {\bm u}, {\bm v},{\bm h}, {\bm s} $
\begin{align}\label{relerr}
\|\bm{w}-\widehat{\bm   w}\|_{rel}=\frac{1}{K}\sum_{k=1}^{K}\frac{\|{\bm  w}^k-\widehat{\bm  w}^k\|_{L^2}}{\|{\bm  w}^k\|_{L^2}}, \quad  \|{\bm  w}^k\|_{L^2}^2=\sum_{i=1}^n\sum_{j=1}^n (w^k_{i,j})^2\Delta x\Delta y,
\end{align}
where $\widehat{\bm  w}=\overline{\bm  w} + V_{w}{\bm w}_r$ denotes the reduced approximation to ${\bm w}$, and $w^k_{i,j}\approx w(x_i,y_j,t_k)$.

Conservation of the discrete conserved quantities \eqref{dener}: the energy, mass, buoyancy,  and total vorticity, of the full and reduced solutions are measured using the time-averaged relative error $\|E\|_{\text{abs}}$ which are given for $E= H,M,B,Q$ and for $\bm{w}=\bm{z},\widehat{\bm{z}}$ as
\begin{equation} \label{conserr1}
\|E\|_{\text{abs}} = \frac{1}{K}\sum_{k=1}^{K} \frac{ |E^k(\bm{w})-E^0(\bm{w})|}{ |E^0(\bm{w})|  }.
\end{equation}

All simulations are performed on a machine with Intel$^{\circledR}$
Core$^{{\mathrm TM}}$ i7 2.5 GHz 64 bit CPU, 8 GB RAM, Windows 10, using 64 bit MatLab R2014.

\begin{figure}[htb!]
	\centering
		\subfloat{\includegraphics[width=0.35\columnwidth]{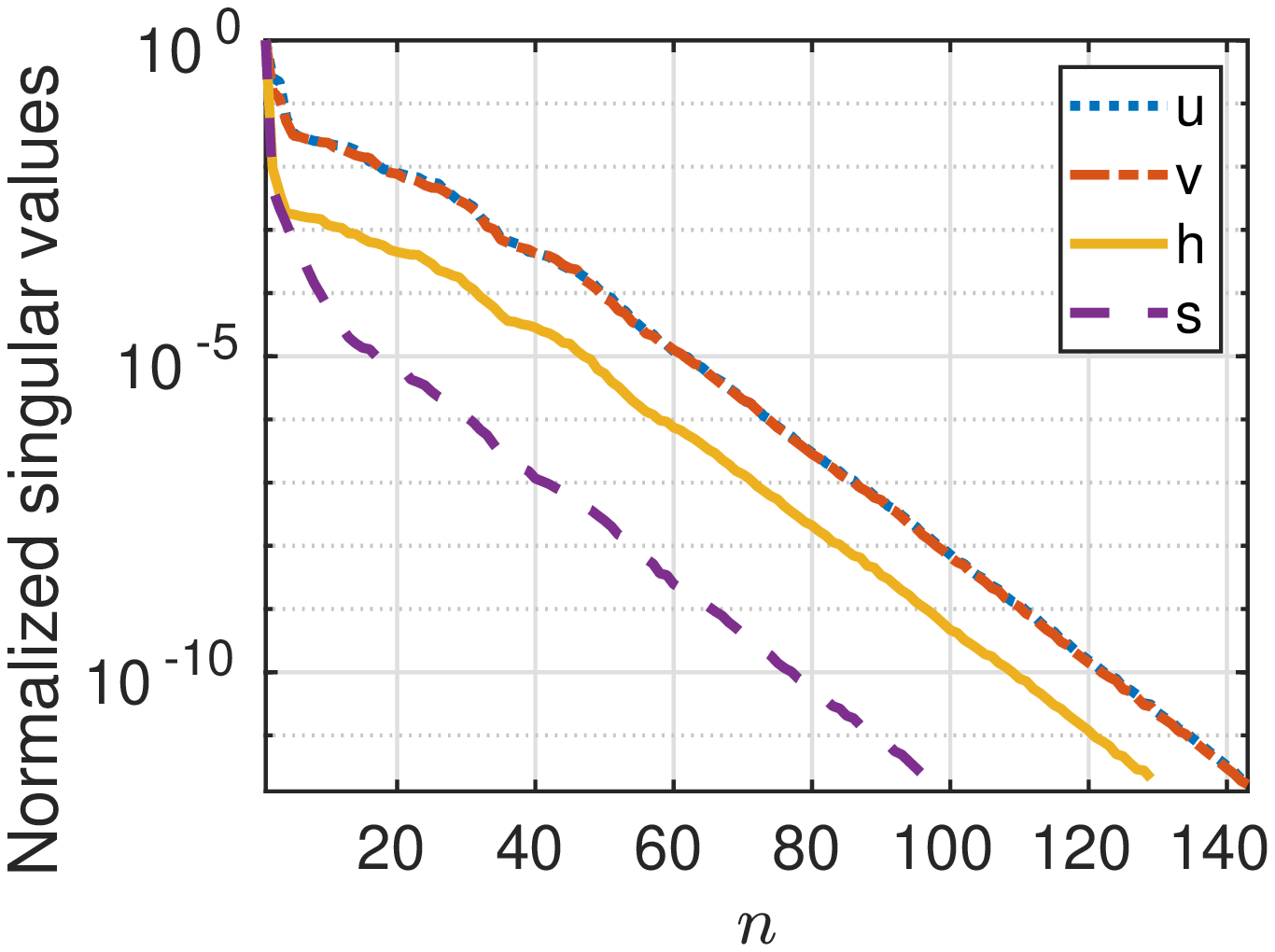}}
		\subfloat{\includegraphics[width=0.35\columnwidth]{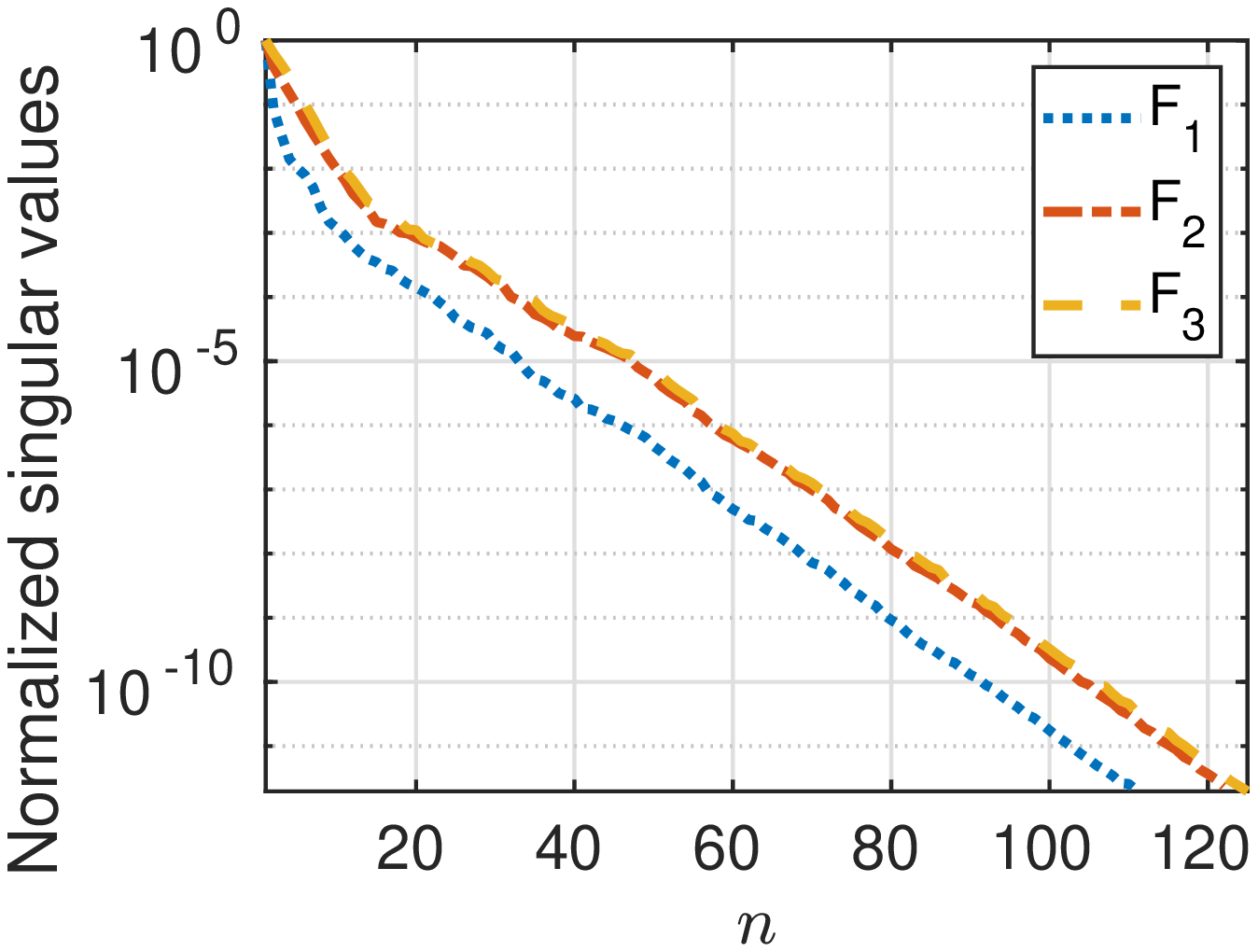}}
	\subfloat{\includegraphics[width=0.35\columnwidth]{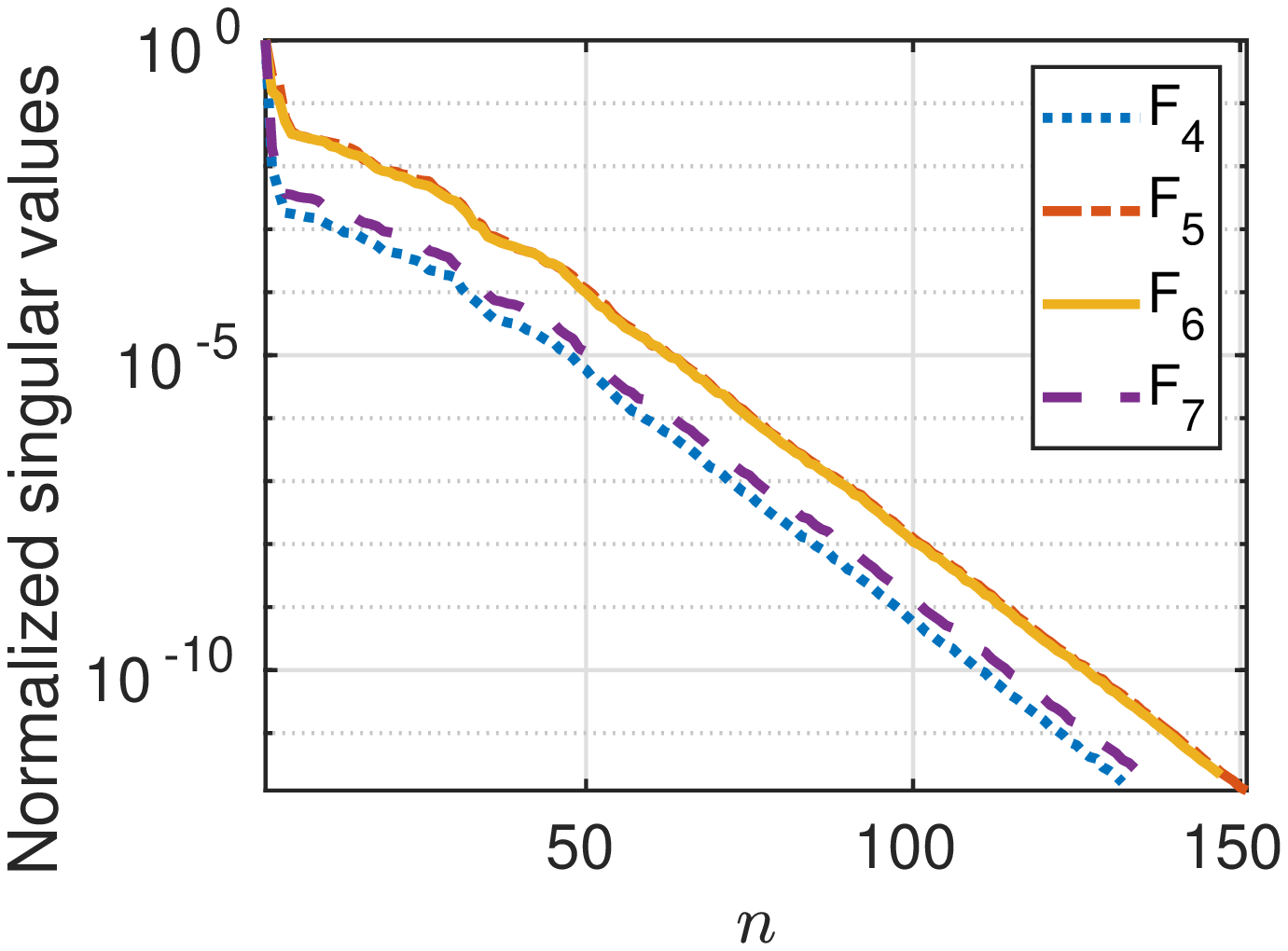}}
		\caption{Normalized singular values: (left) state variables, (middle and right) nonlinear terms  \label{sing}}
\end{figure}

In Figure~\ref{sing}, the singular values decay slowly both for the state variables and  the nonlinear terms, which is the characteristic of the problems with complex wave phenomena in fluid dynamics. According to the energy criteria \eqref{energy_criteria},  $r=5$ POD modes and $p=35$ DEIM modes are selected.

\begin{figure}[htb!]
	\centering
		\subfloat{\includegraphics[width=0.35\columnwidth]{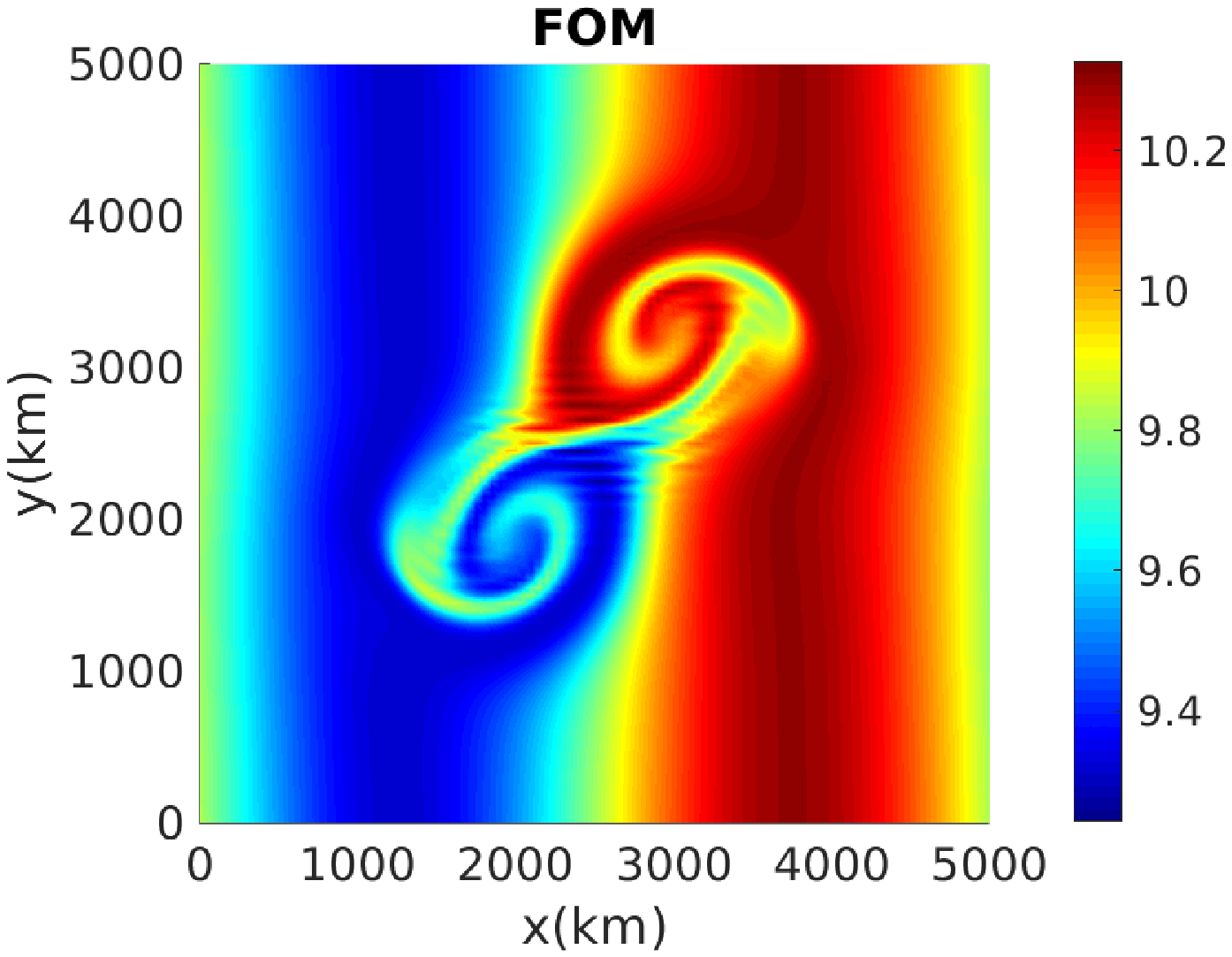}}
		\subfloat{\includegraphics[width=0.35\columnwidth]{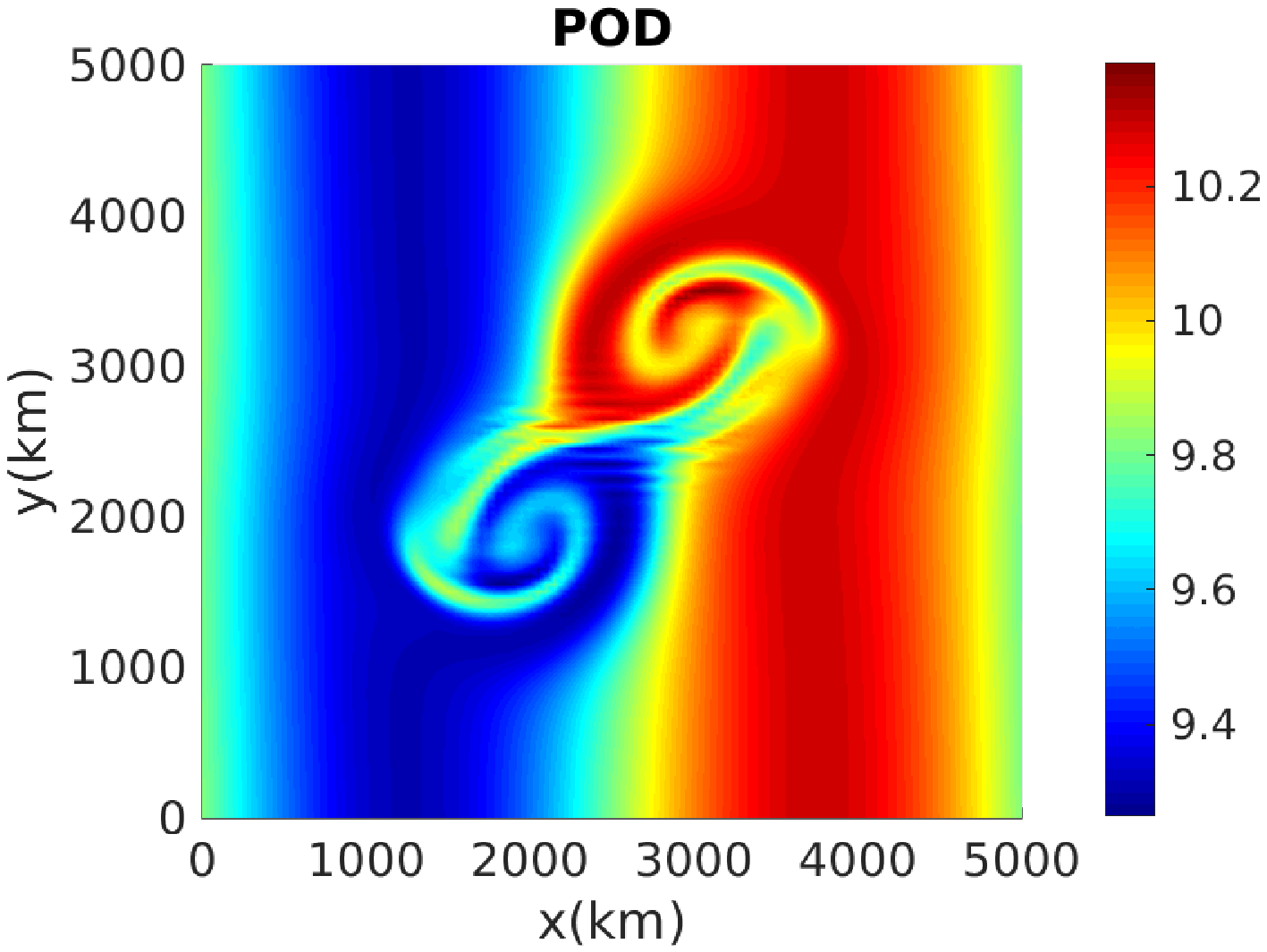}}
		\subfloat{\includegraphics[width=0.35\columnwidth]{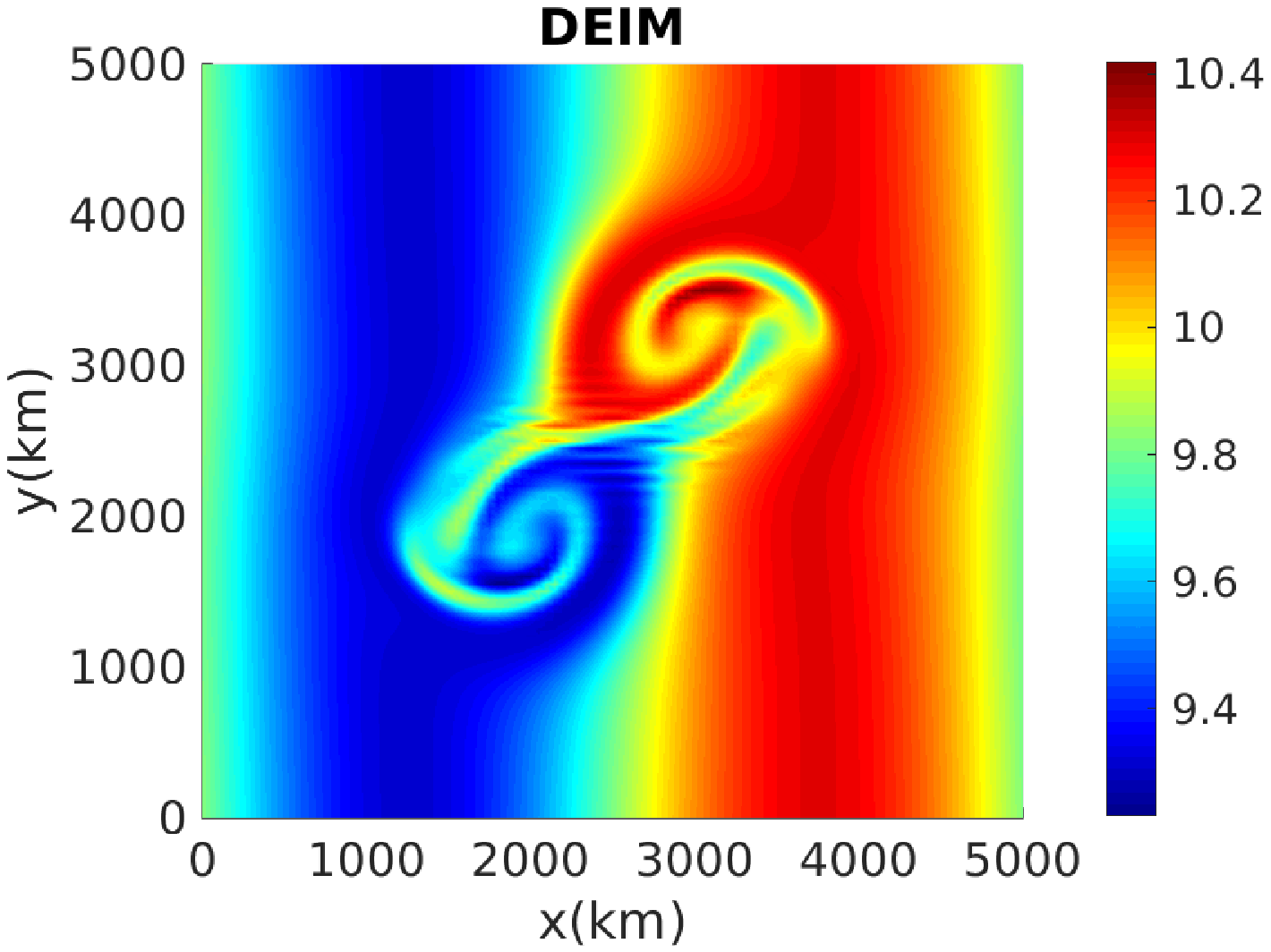}}
		\caption{Bouyancy  at the final time \label{bouyancy}}
\end{figure}

\begin{figure}[htb!]
	\centering
		\subfloat{\includegraphics[width=0.35\columnwidth]{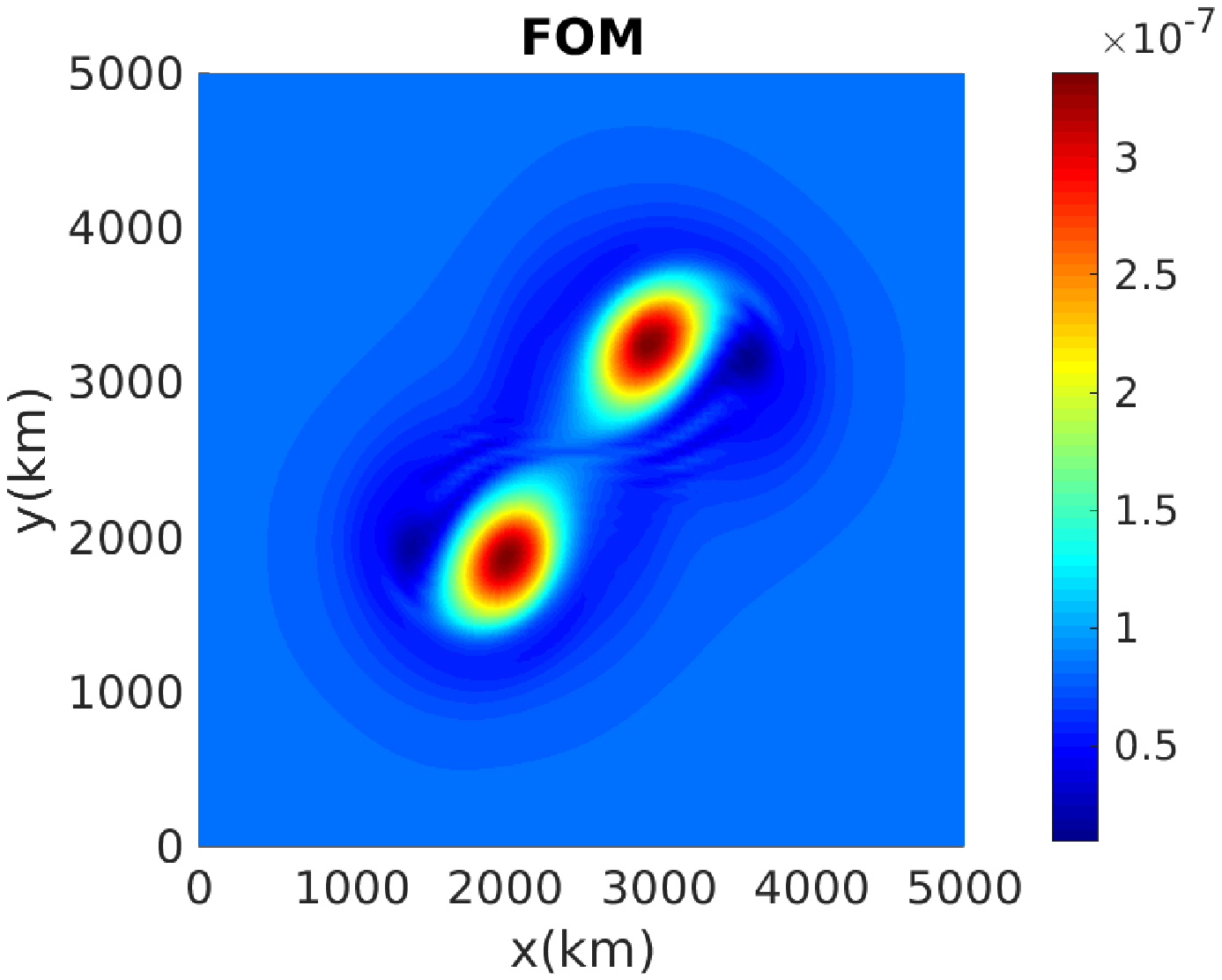}}
		\subfloat{\includegraphics[width=0.35\columnwidth]{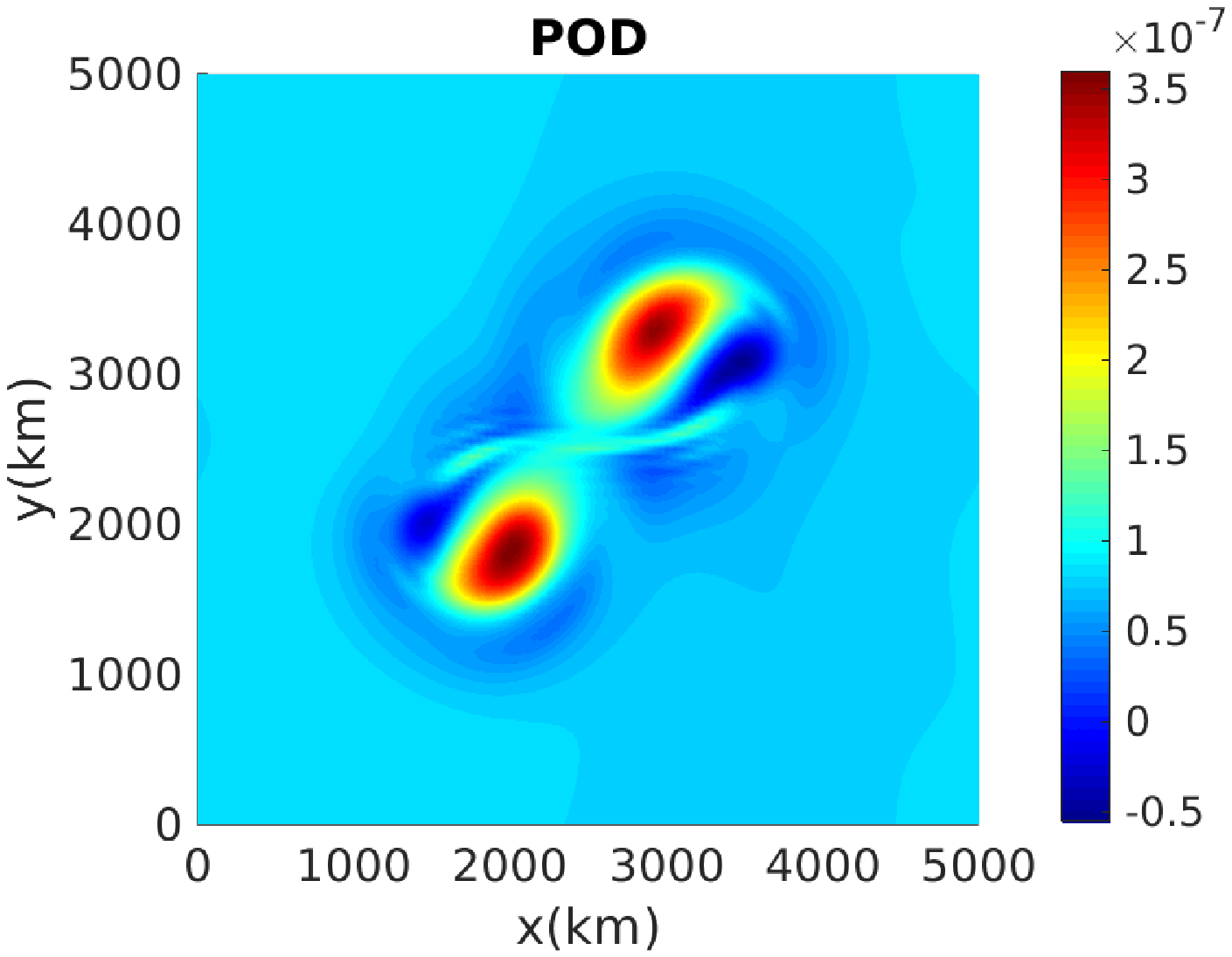}}
		\subfloat{\includegraphics[width=0.35\columnwidth]{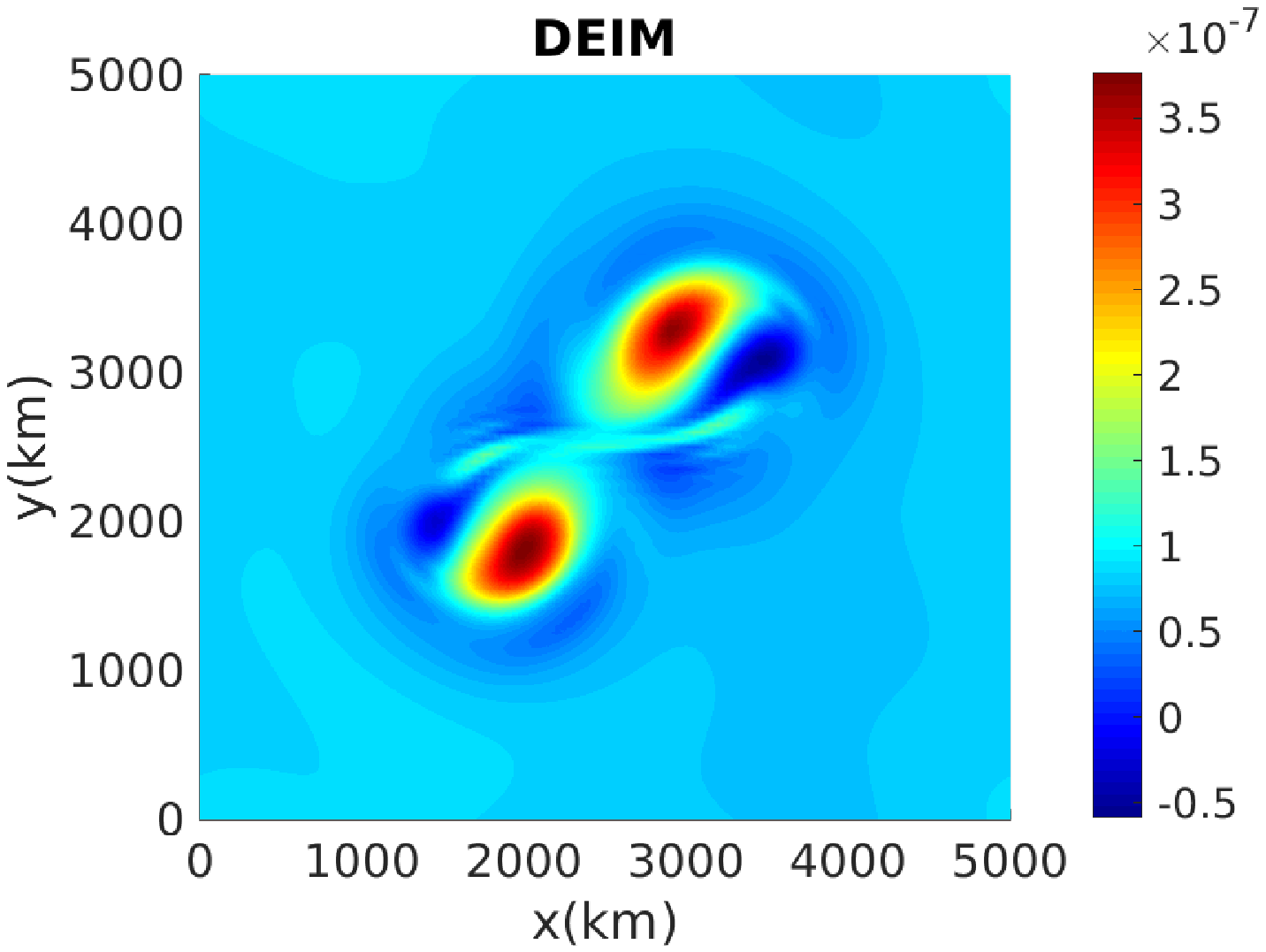}}
		\caption{Potential   vorticity at the final time  \label{vorticity}}
\end{figure}

\begin{figure}[htb!]
	\centering
		\subfloat{\includegraphics[width=0.4\columnwidth]{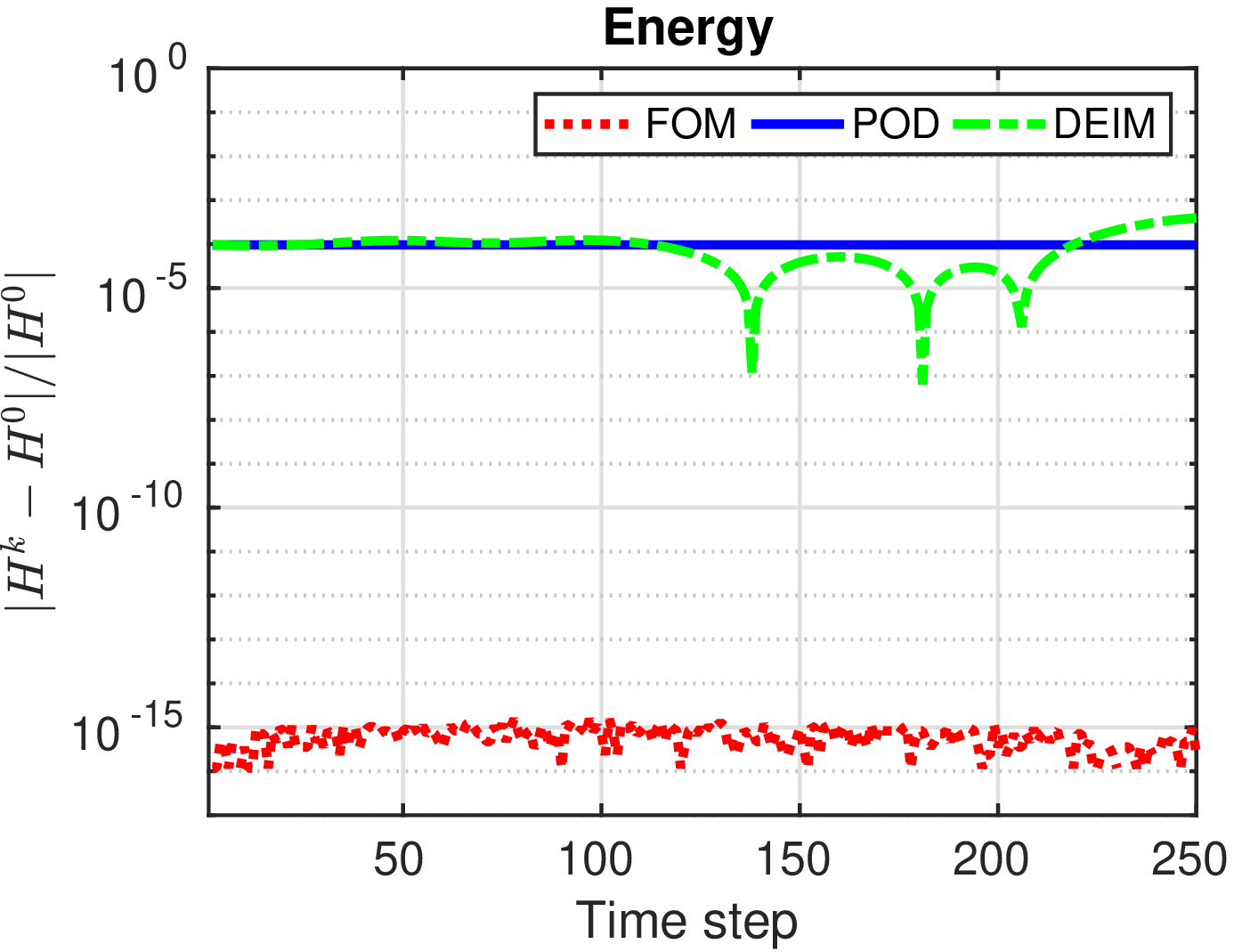}}
		\subfloat{\includegraphics[width=0.4\columnwidth]{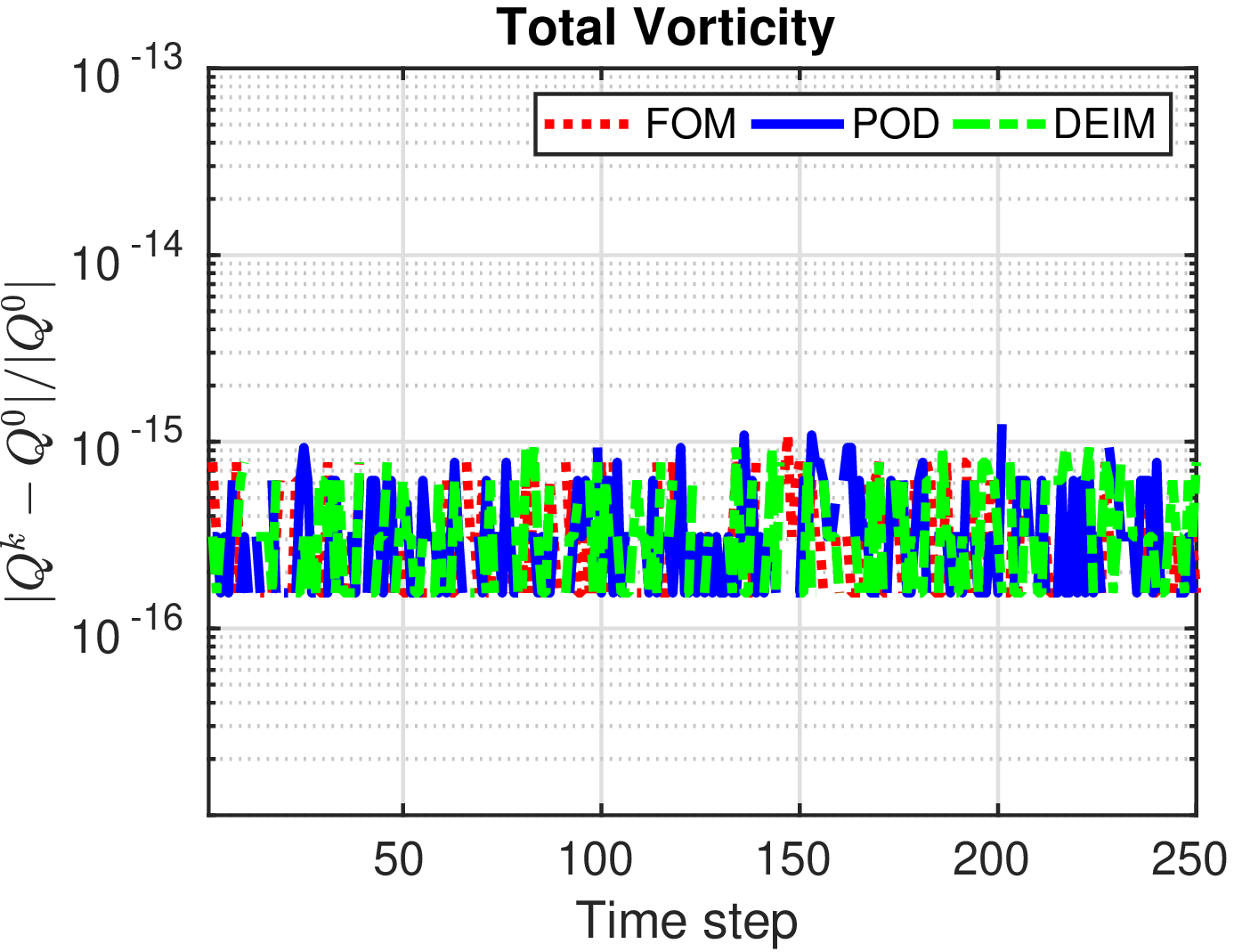}}
		
		\subfloat{\includegraphics[width=0.4\columnwidth]{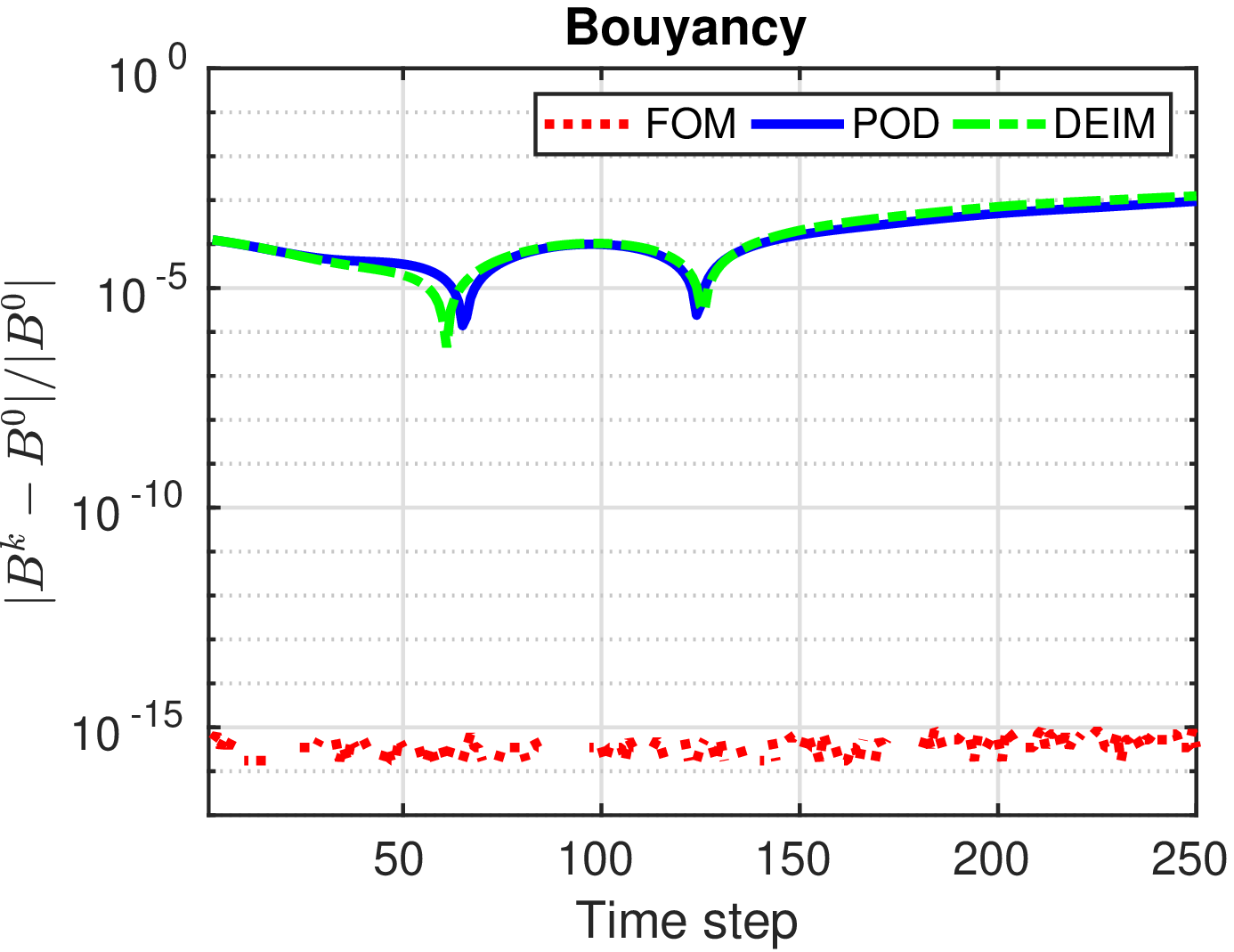}}
		\subfloat{\includegraphics[width=0.4\columnwidth]{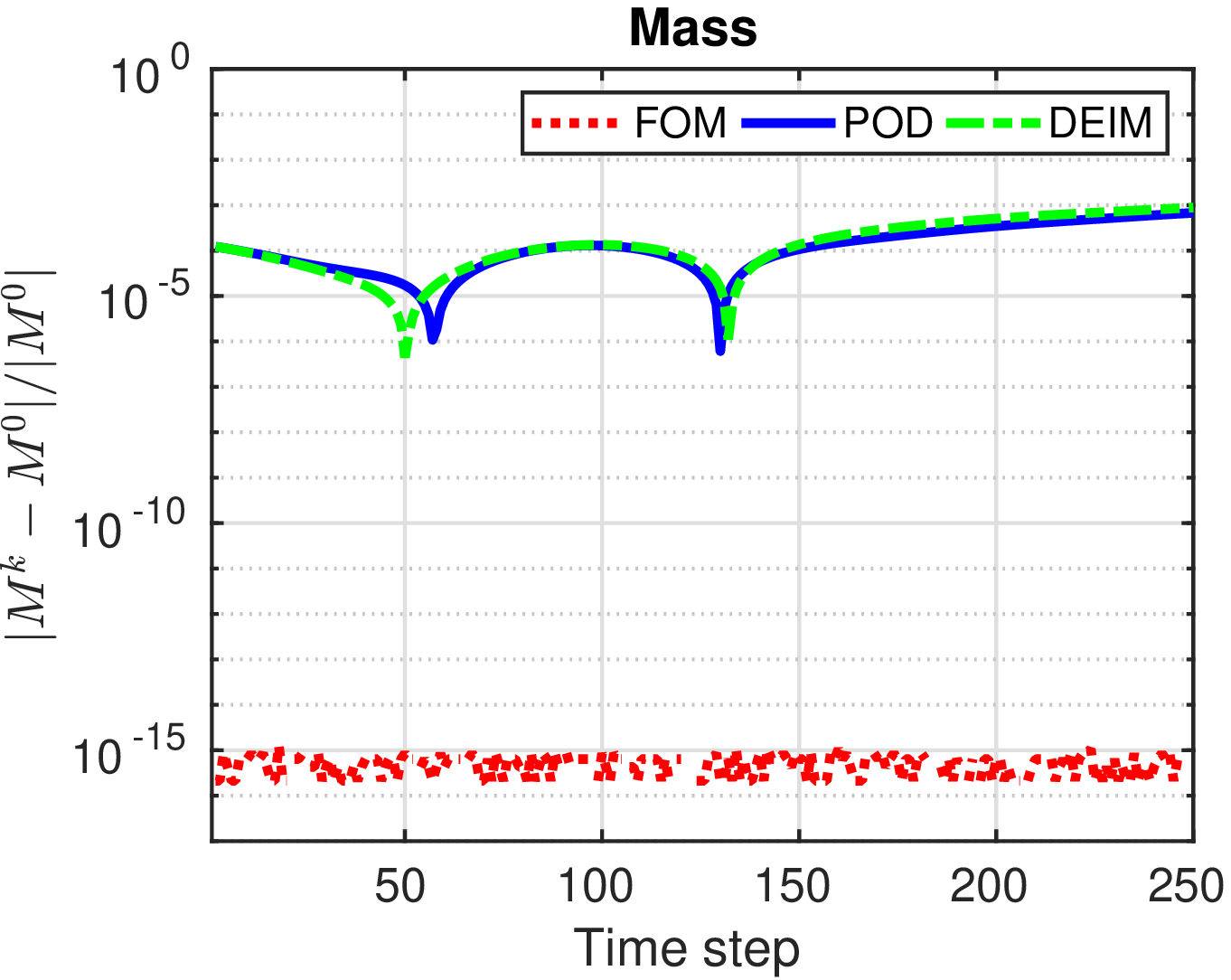}}
		\caption{Relative errors in the conserved quantities. \label{conserved}}
\end{figure}

As shown in the Figures~\ref{bouyancy}-\ref{vorticity}, the bouyancy $\bm{s}$  and the potential vorticity $\bm{q}$ are well approximated by the ROMs at the final time.
In Figure \ref{conserved}, the relative errors in the Hamiltonian (energy)  $|H^k -H^0|/H^0$, the mass  $|M^k - M^0|/M^0$, the buoyancy $|B^k - B^0|/B^0$, and the total  potential vorticity  $|Q^k -Q^0|/Q^0$  versus the time steps  are plotted.
The total potential vorticity is preserved up to machine precision. Energy, mass, and the buoyancy errors of the FOM and ROMs show bounded oscillations over time, i.e., they are preserved approximately at the same
level of accuracy. All conserved quantities are well approximated in the reduced form by the POD and DEIM  and do not show any drift over time. This indicates that the reduced solutions are robust in long term computations.

The time-averaged relative $L_2$ errors of the FOM and ROM solutions in Table~\ref{tbl1} are at the same level of accuracy for
the POD and POD-DEIM. The conserved
quantities are also  accurately preserved by the ROMs as shown  in Table~\ref{tbl2}. The POD-DEIM errors are slightly larger than  the POD errors in
both tables, but the POD-DEIM is much faster than the POD as shown in Table~\ref{tbl3}.

\begin{table}[htb!]
\caption{Time-averaged relative $L_2$ errors of the state variables \label{tbl1}}%
\centering
%\resizebox{\columnwidth}{!}{%
\begin{tabular}{lcccc}
\hline
 &  $\|\bm{h}-\widehat{\bm h}\|_{rel}$ & $\|\bm{u}-\widehat{\bm u}\|_{rel}$  & $\|\bm{v}-\widehat{\bm v}\|_{rel}$  & $\|\bm{s}-\widehat{\bm s}\|_{rel}$ \\
\hline
POD       & 8.622e-03 & 1.502e-01 & 2.185e-01 &  7.899e-04 \\
POD-DEIM  & 1.014e-02 & 1.737e-01 & 2.400e-01 &  7.943e-04  \\
\hline
\end{tabular}
%}
\end{table}

\begin{table}[htb!]
\caption{Mean relative errors of the conserved quantities by FOM and ROMs \label{tbl2}}%
\centering
%\resizebox{\columnwidth}{!}{%
\begin{tabular}{lcccc}
\hline
 & $\|H\|_{\text{abs}}$ & $\|Q\|_{\text{abs}}$  & $\|M\|_{\text{abs}}$ & $\|B\|_{\text{abs}}$ \\
\hline
FOM      & 4.768e-15 &  4.053e-15  &  2.233e-15  &  3.267e-16  \\
POD      & 9.549e-05 &  3.041e-16  &  1.834e-04  &  2.412e-04  \\
POD-DEIM & 9.589e-05 &  3.447e-15  &  2.440e-04  &  3.237e-04 \\
\hline
\end{tabular}
%}
\end{table}

\begin{table}[htb!]
\caption{Wall clock time (in seconds) and speed-up factors \label{tbl3}}
\centering
%\resizebox{\columnwidth}{!}{%
\begin{tabular}{llrr}
\hline
 &  & \textbf{Wall Clock Time} & \textbf{Speed-up} \\
\hline
FOM   &   &   841.2  &   \\
\hline
\multirow{2}{*}{POD}  & offline computation  & 1.7  &  \\
	 & online computation  & 31.5  &   26.7  \\
\hline
 \multirow{2}{*}{POD-DEIM}  & offline computation (POD+DEIM)  & 10.5  &   \\
	 & online computation  &  5.8  & 146.0 \\
\hline
\end{tabular}
\end{table}

In  Table~\ref{tbl3}, we present the computational efficiency results by means of the wall clock time needed for the computations. The time for the offline computations includes the time needed for the basis computation by SVD, computation of the precomputed matrices, and online computation consists of the time needed for the solution of the reduced system and projections to obtain reduced approximations. We note that the wall clock time of the offline computation for POD-DEIM, includes computation of both POD and DEIM basis functions together with the computation of the tensor calculations.
The speed-up factors in Table~\ref{tbl3}, which are calculated as the ratio of the wall clock time needed to compute the FOM over the wall clock time needed for the online computation of the ROM, show that the ROMs with DEIM increases the computational efficiency much further.

%%%%%%%%%%%%%%%%%%%%%%%%%%%%%%%%%%%%%%%%%%%%%%%%%%%%%%%%%%%%%%%%%%%%%%%%%%%%%%%%%%%%%%%
%%%%%%%%%%%%%%%%%%%%%%%%%%%%%%%%%%%%%%%%%%%%%%%%%%%%%%%%%%%%%%%%%%%%%%%%%%%%%%%%%%%%%%%
%%%%%%%%%%%%%%%%%%%%%%%%%%%%%%%%%%%%%%%%%%%%%%%%%%%%%%%%%%%%%%%%%%%%%%%%%%%%%%%%%%%%%%%
\section{Conclusions}
\label{sec:conc}

In this paper, computational efficiency of Hamiltonian \& energy preserving reduced order models are constructed for the RTSWE. We have shown that the  preservation of  the skew-symmetric form  in space and time yields a robust reduced system. Approximating   the Poisson matrix and the gradient of the Hamiltonian with the DEIM, a skew-gradient reduced system is obtained with linear and quadratic terms only.  Application of the tensor techniques to the reduced system, speeds up the computation of the ROMs further.
With relatively small number of POD and DEIM modes, stable, accurate and  fast reduced solutions are obtained;  energy and other conserved quantities are well preserved  over long time-integration. The reduced system can be identified by a reduced energy, and reduced conserved quantities, that mimics those of the high-fidelity system. This results in an, overall, correct evolution of the solutions that ensures robustness of the reduced system.
As a future study, we plan to apply this methodology
to the shallow water magnetohydrodynamic equation \cite{Dellar02,Dellar03}, which  has similar Hamiltonian structure as the RTSWE.

%%%%%%%%%%%%%%%%%%%%%%%%%%%%%%%%%%%%%%%%%%%%%%%%%%%%%%%%%%%%%%%%%%%%%%%%%%%%%%%%%%%%%%%
%%%%%%%%%%%%%%%%%%%%%%%%%%%%%%%%%%%%%%%%%%%%%%%%%%%%%%%%%%%%%%%%%%%%%%%%%%%%%%%%%%%%%%%
%%%%%%%%%%%%%%%%%%%%%%%%%%%%%%%%%%%%%%%%%%%%%%%%%%%%%%%%%%%%%%%%%%%%%%%%%%%%%%%%%%%%%%%
{\bf Acknowledgement}
This work was supported by 100/2000 Ph.D. Scholarship Program of the Turkish Higher Education Council.

%%%%%%%%%%%%%%%%%%%%%%%%%%%%%%%%%%%%%%%%%%%%%%%%%%%%%%%%%%%%%%%%%%%%%%%%%%%%%%%%%%%%%%%
%%%%%%%%%%%%%%%%%%%%%%%%%%%%%%%%%%%%%%%%%%%%%%%%%%%%%%%%%%%%%%%%%%%%%%%%%%%%%%%%%%%%%%%
%%%%%%%%%%%%%%%%%%%%%%%%%%%%%%%%%%%%%%%%%%%%%%%%%%%%%%%%%%%%%%%%%%%%%%%%%%%%%%%%%%%%%%%
\bibliographystyle{unsrtnat}
\bibliography{references}

\end{document}